\documentclass{amsart}
\usepackage{hyperref}
\usepackage{amssymb}
\usepackage{amscd}

\usepackage[english]{babel}

\newcommand{\translation}[1]{$[\![$#1$]\!]$}
\input{tcilatex}
\theoremstyle{definition}
\theoremstyle{remark}
\numberwithin{equation}{section}

\begin{document}
\title[tangency and differentiability]{Tangency vis-\`{a}-vis
differentiability\\
by Peano, Severi and Guareschi}
\author{S.\ Dolecki}
\address{Mathematical Institute of Burgundy\\
CNRS UMR 5584\\
Burgundy University\\
B.P. 47870, 21078 Dijon, France}
\email{dolecki@u-burgogne.fr}
\author{G. H. Greco}
\address{Dipartimento di Matematica, Universit\`{a} di Trento\\
38050 Povo (Tn), Italy}
\email{greco@science.unitn.it}
\dedicatory{Commemorating the 150th Birthday of Giuseppe Peano (1858-1932)}
\date{February 14, 2010.  \emph{To appear in} \textbf{Journal of Convex Analysis}}
\thanks{For the biographical reconstruction related to \textsc{Guareschi} we
are grateful to dott.\,\textsc{Paolo Carrega}, responsible of the Archive
ISRAL, where we could consult the Fondo Guareschi, to dott.\,\textsc{%
Alessandra Baretta} of the Historical Archive of the University of Pavia, to
Ms.\,\textsc{Anna Rapallo} and Ms.\,\textsc{Maddalena De Mola} of the
Historical Archive of the University of Genoa, to Ms. \textsc{Anna Robbiano}
of the CSBMI of the Faculty of Sciences of the University of Genoa, and to
ing. \textsc{Giovanni Paolo Peloso}, the secretary of the Accademia Ligure
delle Scienze e Lettere.}
\maketitle

\begin{abstract}
\textsc{Peano} defined \emph{differentiability} of functions and \emph{lower
tangent cones} in 1887, and \emph{upper tangent cones} in 1903, but uses the
latter concept already in 1887 without giving a formal definition. Both
cones were defined for arbitrary sets, as certain limits of appropriate
homothetic relations. Around 1930\ \textsc{Severi} and \textsc{Guareschi},
in a series of mutually fecundating individual papers, characterized
differentiability in terms of \emph{lower tangent cones} and strict
differentiability in terms of \emph{lower paratangent cones}, a notion
introduced, independently, by \textsc{Severi} and \textsc{Bouligand} in
1928. \textsc{Severi} and \textsc{Guareschi} graduated about 1900 from the
University of Turin, where \textsc{Peano} taught till his demise in 1932.
\end{abstract}

\section{Preamble}

In 2008 mathematical community celebrated the 150th anniversary of the birth of \textsc{%
Giuseppe Peano}, as well as the 100th anniversary of the last (fifth)
edition of \emph{Formulario Mathematico}. Taking part in the commemoration,
we have been reviewing \textsc{Peano}'s foundational contributions to
various branches of mathematics: optimization \cite{DG-Peano}, Grassmann
geometric calculus \cite{greco_grass}, derivation of measures \cite%
{greco-sonia}, definition of surface area \cite{greco-sonia-area}, general
topology \cite{DG-Peano_top}, infinitesimal calculus \cite{GG-libro}, as
well as to tangency and differentiability (in the present paper). \textsc{%
Peano} contributed in an essential way to several other fields of
mathematics: set theory\thinspace\footnote{%
In 1914 \textsc{Hausdorff} wrote in \emph{Grundz\"{u}ge der Mengenlehre} 
\cite[(1914), p.\,369]{hausdorff} of \textsc{Peano}'s filling cruve: \emph{das ist
eine der merkw\"{u}rdigsten Tatsachen der Mengenlehere, deren Entdeckung wir
G. Peano verdanken} \translation{this is one of the most remarkable facts of set theory,
the discovery of which we owe to G.\ Peano}. It is less known that \textsc{%
Peano} formulated the axiom of choice in \cite[(1890)]{peano_equa_diff} (c.f.
Appendix \ref{postcard}), fourteen years before \textsc{Zermelo} \cite[(1904)]%
{zermelo}.}, ordinary differential equations, arithmetic, convexity and,
maybe most significantly, he introduced a completely rigorous formal
language of mathematics. Also these contributions should and hopefully will
be discussed in future papers.

\textsc{Peano} acquired an international reputation soon after his
graduation\thinspace\footnote{%
Already in \cite[(1882)]{peano-super} he observed that the definition of surface measure of the
famous \emph{Cours de calcul diff\'{e}rentiel et int\'{e}gral} of \textsc{%
Serret} \cite{serret} was inadequate.}. Recognized as one
of the leading mathematical authorities of the epoch, he was invited to
publish in prestigious mathematical journals\thinspace\footnote{%
For example, he was invited by \textsc{Klein} to contribute to \emph{%
Mathematische Annalen} (see \textsc{Segre} \cite[(1997)]{segre-michel} and the letters from \textsc{Mayer} to \textsc{Klein} \cite[n.\,125 p.\,161, n.\,126 p.\,163, n.\,148 p.\,181]{mayer}). As a result, \textsc{Peano}
published three papers: on the resolvent (in particular, the exponential of
a matrix) of a system of linear differential equations \cite[(1888)]%
{peano_equa_diff_lin}, on the existence of solutions of a system of
differential equations with the sole hypothesis of continuity \cite[(1890)]%
{peano_equa_diff}, and on a filling curve \cite[(1890)]{peano_curva}.}. He was at
the summit of fame at the break of the 20th century when he took part in the
International Congress of Philosophy and the International Congress of
Mathematicians in Paris in 1900. \textsc{Bertrand Russell}, who also
participated in the philosophy congress, noted in \cite[(1967), pp.\,217-218]{BRm}

\begin{quotation}
The Congress was a turning point in my intellectual life, because I there
met Peano. [...] In discussions at the Congress I observed that he was
always more precise than anyone else, and that he invariably got the better
of any argument upon which he embarked.
\end{quotation}

In \emph{The Principles of Mathematics} \cite[(1903) p.\,241]{BR} \textsc{%
Russell} said that \textsc{Peano} had a rare immunity from error.

\textsc{Peano} was associated with the University of Turin during his whole
mathematical career, from October 1876, when he became a student, till 19th
of April 1932, when he taught his classes as usual, a day before his death.
From 1903 on, following the example of \textsc{M\'{e}ray}, with whom he
corresponded, \textsc{Peano} dedicates himself more and more to auxiliary
international languages (postulated as \emph{lingua rationalis} by \textsc{%
Leibniz} \cite[(1901), Ch.\,III]{couturat}) in company with a philosopher and
logician \textsc{Louis Couturat}, linguists \textsc{Otto Jespersen} and 
\textsc{Jan Baudouin de Courtenay}, and a chemist \textsc{Wilhem Ostwald}\thinspace\footnote{\textsc{Wilhem Ostwald} (1853-1932), Nobel Prize in Chemistry in
1909, in his \emph{Selbstbiographie} \cite[(1927)]{ostwald} describes \textsc{Peano}
as follows:
\par
\begin{quotation}
Eine Personalit\"{a}t besonderer Art war der italienische Mathematiker
Peano. Lang, \"{a}u\ss erst mager, nach Haltung und Kleidung ein
Stubengelehrter, der f\"{u}r Nebendinge keine Zeit hat, mit gelbbleichem,
hohlem Gesicht und tiefschwarzem, sp\"{u}rlichem Haar und Bart, erschien er
ebenso abstrakt, wie seine Wissenschaft. Er hatte eigene Vorschl\"{a}ge zu
vertreten, n\"{a}mlich sein latino sine flexione, ein tunlichst
vereinfachtes latein, f\"{u}r welches er mit unersch\"{u}tterlicher Hingabe
eintrat, da er als Italiener das Gef\"{u}hl hatte, im Latein ein uraltes
Erbe zu verteidigen.
\end{quotation}

\translation{An Italian mathematician Peano was a personality of peculiar kind. Tall, extremely slim, by attitude and clothes, a scientist, who has no time for secondary things, with his pale yellowish hollow face and sparse deeply black hair and beard, looked so abstract as his science. He had his proper proposal to present, namely his latino sine flexione, a simplified, as much as possible, Latin, which he presented with imperturbable devotion, since, as an Italian, he had the feeling to defend in Latin a primordial heritage.}}.
 This interest becomes his principal passion after the
completion of the last edition of \emph{Formulario Mathematico} in 1908,
written in a (totally rigorous) mathematical formal language\thinspace\footnote{%
\textsc{Hilbert} and \textsc{Ackermann} write in the introduction to \cite[(1928)]%
{hilbert-Ackermann}: \emph{G. Peano and his co-workers began in 1894 the
publication of the Formulaire de Math\'{e}matiques, in which all the
mathematical disciplines were to be presented in terms of the logical
calculus.}} and commented in an auxiliary language, \emph{latino sine
flexione}, both conceived by \textsc{Peano}.

It should be emphasized that the formal language conceived and used by 
\textsc{Peano} was not a kind of shorthand adapted for a mathematical
discourse, but a collection of ideographic symbols and syntactic rules with
univocal semantic interpretations, which produced precise mathematical
propositions, as well as inferential rules that ensure the correctness of
arguments.

\textsc{Peano}'s fundamental contributions to mathematics are numerous. Yet,
nowadays, only few mathematical achievements are commonly associated with
his name. It is dutiful to reconstitute from (partial) oblivion his
exceptional role in the development of science (see Appendix \ref{postcard}%
). In the present paper we intend to delineate the evolution, in the work of 
\textsc{Peano}, of the concept of tangency and of its relation to
differentiability\thinspace\footnote{%
In his reference book \cite[(1973)]{May} \textsc{K.O.\thinspace May} discusses a
role of direct and indirect sources in historiography of mathematics. He
stresses the importance of primary sources, but acknowledges also the
usefulness of secondary (and $n$-ary sources) under the provision of
critical evaluation. As mathematicians, we are principally interested in
development of mathematical ideas, so that we use almost exclusively primary
sources, that is, original mathematical papers. On the other hand, one
should not neglect the biography of the mathematicians whose work one
studies, because it provides information about effective and possible
interactions between them.}.

By respect for historical sources and for the reader's convenience, the
quotations in the sequel will appear in the original tongue with a
translation in square brackets (usually placed in a footnote). All the
biographical facts concerning \textsc{Peano} are taken from \textsc{%
H.C.~Kennedy}, \emph{Life and Works of Giuseppe Peano} \cite[(1980, 2006)]%
{kennedy,kennedy_corrected}. On the other hand, we have checked all the reported
bibliographic details concerning mathematical aspects.

\section{Introduction}

In \emph{Applicazioni Geometriche} of 1887 \cite{peano87}, \textsc{Peano}
defined \emph{differentiability} of functions, \emph{lower tangent cone},
and (implicitly in \cite{peano87} and explicitly in \emph{Formulario
Mathematico} of 1903 \cite{formulaire4}) \emph{upper tangent cone}, both for
arbitrary sets, as certain limits of appropriate homothetic relations.
Around 1930\ \textsc{Francesco Severi} (1879-1961) and \textsc{Giacinto
Guareschi} (1882-1976), in a series of mutually fecundating individual
papers, characterized differentiability in terms of tangency without
referring to \textsc{Peano}.

Following \textsc{Peano} \cite[(1908) p.\,330]{peano_1908}, a function $%
f:A\rightarrow \mathbb{R}^{n}$ is \emph{differentiable at} an accumulation
point\emph{\ }$\hat{x}$ of $A\subset \mathbb{R}^{m}$ if $\hat x\in A$ and there exists\thinspace\footnote{%
In his definition \textsc{Peano} assumes uniqueness, which we drop because
of the prevalent contemporary use that we adopt in the sequel of the paper.}
a linear function $Df(\hat{x}):\mathbb{R}^{m}\rightarrow \mathbb{R}^{n}$
such that%
\begin{equation}
\lim\nolimits_{A\ni x\rightarrow \hat{x}}\dfrac{f(x)-f(\hat{x})-Df(\hat{x}%
)(x-\hat{x})}{\Vert x-\hat{x}\Vert }=0.  \label{differ}
\end{equation}%
It is \emph{strictly differentiable at }$\hat{x}$ (\textsc{Peano} \cite[%
(1892)]{peano_mathesis} for $n=1$, \textsc{Severi} in \cite[(1934) p.\,185]%
{Severi-diff}\thinspace\footnote{%
As we will see later, \textsc{Severi} uses the term \emph{hyperdifferentiable%
}.}) if (\ref{differ}) is strengthened to

\begin{equation}
\lim\nolimits_{A\ni x,y\rightarrow \hat{x},{x\neq y}}\dfrac{f(y)-f(x)-Df(%
\hat{x})(y-x)}{\Vert y-x\Vert }=0.  \label{strict differ}
\end{equation}%
These are exactly the definitions that we use nowadays. The first notion is
often called \emph{Fr\'{e}chet differentiability }(referring to \textsc{%
Fr\'echet} \cite[(1911)]{Frechet_diff,Frechet_diff2}) and the second is
frequently referred to \textsc{Leach} \cite[(1961)]{leach}, where it is
called \emph{strong differentiability}
and to  \textsc{Bourbaki} \cite[(1967), p.\,12]{bourbaki}.

Currently an assortment of tangent cones have been defined by a variety of
limits applied to homothetic relations. \textsc{Peano} gave an accomplished
definition of tangency in \emph{Formulario Mathematico}\ \cite[(1908)]%
{peano_1908}, as was noticed in \textsc{Dolecki}, \textsc{Greco} \cite[(2007)]%
{DG-Peano}; he defined what we call respectively, the \emph{lower} and the 
\emph{upper tangent cones} of $F$ at $x$ (traditionally denominated \emph{%
adjacent} and \emph{contingent} cones)\thinspace\footnote{%
Actually Peano defined affine variants of these cones.}%
\begin{eqnarray}
\limfunc{Tan}\nolimits^{-}(F,x):= \limfunc{Li}\limits_{t\rightarrow 0^{+}}%
\tfrac{1}{t}\left( F-x\right) ,  \label{lower} \\
\limfunc{Tan}\nolimits^{+}(F,x):= \limfunc{Ls}_{t\rightarrow 0^{+}}\tfrac{1%
}{t}\left( F-x\right) ,  \label{upper}
\end{eqnarray}%
where $\limfunc{Li}\limits_{t\rightarrow 0^{+}}$ and $\limfunc{Ls}%
\limits_{t\rightarrow 0^{+}}$ denote the usual lower and upper limits of
set-valued maps. Here, we adopt the modern
definition of lower and upper limits in metric spaces, both introduced by 
\textsc{Peano}, the first in \emph{Applicazioni geometriche} \cite[(1887),  p.\,302]{peano87} and the second 
in \emph{Lezioni di analisi infinitesimale} \cite[(1893), volume 2, p.\,187]{peano1893}
(see \textsc{Dolecki}, \textsc{Greco} \cite[(2007)]{DG-Peano} for
further details). Let $d$ denote the Euclidean distance on $\mathbb{R}^{n}$
and let $A_{t}$  be a subset of $\mathbb{R}^{n}$ for $t>0$.
According to 
\textsc{Peano}, 
\begin{eqnarray}
 \limfunc{Li}\nolimits_{t\rightarrow 0^{+}}A_{t}:=\{x\in 
\mathbb{R}^{n}:\lim_{t\rightarrow 0^{+}}d(x,A_{t})=0\}  \label{lowerlimit} \\
\limfunc{Ls}%
\nolimits_{t\rightarrow 0^{+}}A_{t}:=\{x\in \mathbb{R}^{n}:\liminf_{t%
\rightarrow 0^{+}}d(x,A_{t})=0\}. \label{upperlimit}
\end{eqnarray}%
Since $d(v,\frac{1}{t}(F-x))=\frac{1}{t}%
d(x+tv,F)$, from   \eqref{lower} and \eqref{upper} it follows that 
\begin{eqnarray}
v\in \limfunc{Tan}\nolimits^{-}(F,x)\quad \text{if and
only if}\quad\lim_{t\rightarrow 0^{+}}\frac{1}{t}d(x+tv,F)=0\\
v\in \limfunc{%
Tan}\nolimits^{+}(F,x)\quad \text{if and
only if}\quad\liminf_{t\rightarrow 0^{+}}\frac{1}{%
t}d(x+tv,F)=0.
\end{eqnarray}%

The \emph{upper paratangent cone }(traditionally called \emph{paratingent
cone}) of $F$ at $x$%
\begin{equation}
\limfunc{pTan}\nolimits^{+}(F,x):=\limfunc{Ls}\limits_{t\rightarrow
0^{+},\;F\ni y\rightarrow x}\tfrac{1}{t}\left( F-y\right)  \label{para}
\end{equation}%
was introduced later by \textsc{Severi} \cite[(1928) p.\,149]{Severi-conf-RM}
and \textsc{Bouligand} in \cite[(1928) pp.\,29-30]{boul_BSMF}\thinspace
\footnote{%
Successively in \cite[(1930) pp.\,42-43]{boul_RS}, \textsc{Bouligand}
introduces the terms of \emph{contingent} and \emph{paratingent} to denote
upper tangent and paratangent cones. In contrast to definitons \eqref{upper}
and \eqref{para}, for \textsc{Severi} and \textsc{Bouligand}, an upper
tangent (resp. upper paratangent) cone is a family of half-lines (resp.
straight lines); consequently, they are empty at isolated points and, on the
other hand, they consider closedness in the sense of half-lines (resp. straight lines).}. The \emph{%
lower paratangent cone}%
\begin{equation}
\limfunc{pTan}\nolimits^{-}(F,x):=\limfunc{Li}\limits_{t\rightarrow
0^{+},\;F\ni y\rightarrow x}\tfrac{1}{t}\left( F-y\right)  \label{clarke}
\end{equation}%
is usually called the \emph{Clarke tangent cone } (see \textsc{Clarke} \cite[(1973)]%
{clarke}). In \cite[pp.\,499-500]{DG-Peano} we listed the properties of the
upper tangent cone observed by \textsc{Peano}. Of course, \qquad 
\begin{equation}
\limfunc{pTan}\nolimits^{-}(F,x)\subset \limfunc{Tan}\nolimits^{-}(F,x)%
\subset \limfunc{Tan}\nolimits^{+}(F,x)\subset \limfunc{pTan}%
\nolimits^{+}(F,x).  \label{4cones}
\end{equation}

In the works of \textsc{Peano} there are no occurrences of sets for which
the upper and lower tangent cones are different. Here we furnish an easy 
one.\thinspace\footnote{%
\label{footTan}In \cite[(2007), p.\,499, footnote 21]{DG-Peano} we observed that $v\in \limfunc{Tan}%
\nolimits^{-}(S,x)$ if and only if 
\begin{itemize}
\item[(*)] there exists a sequence $\left\{
x_{n}\right\} _{n}\subset S$ such that $\lim\nolimits_{n}x_{n}=x$ and $%
\lim\nolimits_{n}n(x_{n}-x)=v$. 
\end{itemize}
On the other hand it is well known that $v\in \limfunc{Tan}%
\nolimits^{+}(S,x)$ if and only if 

\begin{itemize}
\item[(**)] there exist sequences $\{\lambda_{n}\}_{n}\subset\mathbb{R}_{++}$ and $\left\{
x_{n}\right\} _{n}\subset S$ such that  $\lim\nolimits_{n}\lambda_{n}=0$ , $\lim\nolimits_{n}x_{n}=x$ and $%
\lim\nolimits_{n}(x_{n}-x)/\lambda_{n}=v$. 
\end{itemize}
In \cite[(1929)]{vonN} \textsc{von Neumann} shows that a closed matrix group $G$  is a Lie group whenever $(1)$  $\limfunc{Tan}\nolimits^{+}(G,E)$ at the unit $E$ of $G$  is a matrix Lie algebra, $(2)$ $(**)$ implies $(*)$ and $(3)$  $\exp A\in G$ for every $A\in \limfunc{Tan}\nolimits^{-}(G,E)$.  The second claim, which amounts to $\limfunc{Tan}\nolimits^{+}(G,E)=\limfunc{Tan}\nolimits^{-}(G,E)$, is the crucial step in his proof.}

\begin{example}
If $S:=\left\{ \frac{1}{n!}:n\in \mathbb{N}\right\} $, then $\limfunc{Tan}%
\nolimits^{+}(S,0)=\mathbb{R}_{+}$ and $\limfunc{Tan}\nolimits^{-}(S,0)=%
\left\{ 0\right\} $.
\end{example}

It is surprising, but it seems that so far in the literature there have been
no such examples. The pretended instances: 
\begin{equation*}
A:=\{(t,t\sin (\frac{1}{t})):t\in \mathbb{R}\smallsetminus \{0\}\}
\end{equation*}%
given by \textsc{Rockafellar} and \textsc{Wets} in \emph{Variational Analysis%
} \cite[(1998), p.\,199]{rock-wets}, and 
\begin{equation*}
B:=\{(t,-t):t<0\}\cup \{(\frac{1}{n},\frac{1}{n}):n\in \mathbb{N}\}
\end{equation*}%
provided by \textsc{Aubin} and \textsc{Frankowska} in \emph{Set-Valued
Analysis} \cite[(1990), p.\,161]{aubin-frank} are not pertinent, because in both of
them the upper and the lower tangent cones coincide\thinspace\footnote{%
In fact, by footnote \ref{footTan}, $\limfunc{Tan}\nolimits^{+}(A,(0,0))=%
\limfunc{Tan}\nolimits^{-}(A,(0,0))=\left\{ (h,k)\in \mathbb{R}%
^{2}:\left\vert k\right\vert \leq \left\vert h\right\vert \right\} $ and $%
\limfunc{Tan}\nolimits^{+}(B,(0,0))=\limfunc{Tan}\nolimits^{-}(B,(0,0))=%
\left\{ (t,\left\vert t\right\vert ):t\in \mathbb{R}\right\} $.}.

In the literature there are numerous examples of sets, for which other
inclusions in (\ref{4cones}) are strict.

The remarkable fact that the coincidence of the upper and lower paratangent cones at every point of a
locally closed subset $F$ of Euclidean space is equivalent to the fact
that $F$ is a $C^{1}$-submanifold, has not been observed till now.
It will be an object of \cite{manifold}, in which a  mathematical and historical account on the subject will be provided.
\thinspace\footnote{ Although  \textsc{Severi} and \textsc{Guareschi} characterized $C^{1}$ manifolds in Euclidean space in terms of tangency, their definitions and reasonings are not entirely transparent; see \textsc{Greco} \cite{manifold} for further details.}

Intrinsic notions of tangent straight line to a curve and of tangent plane
to a surface were clear to \textsc{Peano} (see Section \ref{evo}) and even
prior to him, before the emergence of the concept of tangent cone to an
arbitrary set. On rephrasing these special notions in terms of a vector
space $H$, tangent to a set $F$ at an accumulation point $\hat{x}$ of $F$, we
recover the following condition:%
\begin{equation}
\lim\nolimits_{F\ni x\rightarrow \hat{x},x\neq \hat{x}}\dfrac{d(x,H+\hat x)}{%
d\left( x,\hat{x}\right) }=0.  \label{tangtrad}
\end{equation}%
Geometrically, (\ref{tangtrad}) means that the vector space $H$ and the half-line passing through $%
\hat{x}$ and $x$ in $F$  form an angle that tends to zero as $x$ tends $\hat{x%
}$.

From 1880 \textsc{Peano} taught at the University of Turin.
Among the students of that university at the very end of 19th century were 
\textsc{Beppo Levi}, \textsc{Severi} and \textsc{Guareschi} (see the
biography in Appendix \ref{biographies}). They were certainly acquainted
with the famous \emph{Applicazioni Geometriche} \cite[(1887)]{peano87} of 
\textsc{Peano}, so that their writings on tangency and differentiability
could not abstract from the achievements of \textsc{Peano}. Yet neither 
\textsc{Severi} nor \textsc{Guareschi} cite \textsc{Peano}\thinspace\footnote{\label%
{ftnote}\textsc{Severi} however mentions in \cite[(1930)]{Severi-Krak} a
paper \cite[(1930)]{Cassina-Milano} of \textsc{Cassina} (who, by the way, became
later the editor of the collected works of \textsc{Peano} \cite{opere}). It
was on browsing through \textsc{Severi}'s citation of \textsc{Cassina} that
the second author (\textsc{G.\,H.\,Greco}) of this paper discovered the immensity of \textsc{Peano}'s
contributions to scientific culture. Parenthetically, \textsc{Severi}
reproaches to \textsc{Cassina} for having failed to quote him:
\par
\begin{quotation}
[\textsc{Cassina}] ha ultimamente considerato allo stesso mio modo la figura
tangente ad un insieme, ignorando certo i precedenti sull'argomento.
\end{quotation}
\par
\translation{[\textsc{Cassina}] recently considered, in the same way of mine, the
tangent figure of a set, apparently ignoring the precedents in this topic.}
\par
This surprising oblivion of \textsc{Peano}'s work by \textsc{Severi} can be
perhaps explained by a merely sporadic interest in mathematical analysis by
this algebraic geometer.
\par
Another algebraic geometer, \textsc{Beniamino Segre} (a coauthor with 
\textsc{Severi} of a paper on tangency \cite[(1929)]{Severi-para}, and, on
the other hand, an author of a historical paper on \textsc{Peano} \cite[%
(1955)]{segre}), presented to \emph{Accademia dei Lincei} a paper on
tangency \cite[(1973)]{ursescu} that ignored the contributions of \textsc{%
Peano}, \textsc{Severi} and \textsc{Segre} himself, without reacting to this
unawareness.
\par
It is also surprising that \textsc{Boggio}, one of the best known pupils of 
\textsc{Peano}, did not recall in \cite[(1936)]{boggio} the famous
contribution to tangency of his mentor, when he recommended for publication
in \emph{Memorie dell'Accademia delle Scienze di Torino} a paper of \textsc{%
Guareschi} \cite[(1936)]{Guareschi-diff-To} that begins: ``\emph{Il concetto
di semitangente [...] introdotto nell'analisi da F. Severi.}'' \translation{The concept of
semitangent [...] introduced in analysis by F. Severi.}}.
By the bye, in \cite%
[(1932)]{LeviB} \textsc{Beppo Levi} acknowledges explicitly the influence of 
\emph{Calcolo Geometrico} \cite[(1888)]{calcolo88} of \textsc{Peano} on his
understanding of the work of \textsc{Grassmann}; \textsc{Beppo Levi} recalls
his enthusiastic interest in \emph{Calcolo Geometrico} and difficulty in
reading \emph{Ausdehnungslehre} \cite{grass}:

\begin{quotation}
[interesse] quasi entusiastico che, giovane principiante, mi prese alla
lettura del \emph{Calcolo geometrico secondo l'Ausdehnungslehre di Grassmann}%
; e ricordo all'opposto, l'impressione di malsicura astrattezza che il
medesimo principiante ricevette volendo affrontare la fonte, l'\emph{%
Ausdehnungslehre} del 1844.\thinspace\footnote{%
\translation{Almost entusiatic [interest] that took me, a young beginner, at the lecture
of \emph{Calcolo geometrico secondo l'Ausdehnungslehre di Grassmann}; and I
remember, in contrast, an impression of insecure abstractness that the same
beginner received attempting to confront the source, \emph{Ausdehnungslehre}
of 1844.}}
\end{quotation}

In \cite[p.\,241]{Frechet_3} of 1937, \textsc{Fr\'{e}chet} comments\thinspace\footnote{%
We believe that \textsc{Fr\'{e}chet}, who never investigated tangency, took
this information either from his friend \textsc{Bouligand} or, directly,
from a paper of \textsc{Severi} \cite[(1931)]{Severi-Krak} where \textsc{B.\,Levi}, \textsc{Bouligand} and his pupils \textsc{Rabat\'e} and \textsc{Durand} are quoted. To our knowledge \textsc{%
Bouligand} neither refers to nor quotes \textsc{Severi}.}:

\begin{quotation}
On doit \`{a} M. Bouligand et \`{a} ses \'{e}l\`{e}ves d'avoir entrepris l'%
\'{e}tude syst\'{e}matique [de] cette th\'{e}orie des \textquotedblleft
contingents et paratingents\textquotedblright\ dont l'utilit\'{e} a \'{e}t%
\'{e} signal\'{e}e d'abord par M. Beppo Levi, puis par M. Severi.\thinspace\footnote{%
\translation{We owe to Bouligand and his pupils a systematic study [of] this theory of
contingents and paratingents, the usefulness of which was pointed out first
by Beppo Levi, then by Severi.}}
\end{quotation}

Following the guidelines of \textsc{Fr\'{e}chet}, we initiated to study the
writings of \textsc{Severi} (see, for example, \textsc{Dolecki} \cite[(1982)]{tang-diff}) and, thanks
to a reference in \textsc{Severi} \cite[(1934)]{Severi-diff}, also those of \textsc{%
Guareschi}.

An exhaustive historical study of the work of \textsc{Bouligand} and his
pupils is also dutiful, and we hope that it will be done before long\thinspace\footnote{%
Among those who refer to \textsc{Bouligand} in their study of tangency we recall \textsc{Durand},
\textsc{Rabat\'{e}} (1931), \textsc{Mirguet} (1932), \textsc{Marchaud}
(1933), \textsc{Blanc} (1933), \textsc{Charpentier} (1933), \textsc{Vergn%
\`{e}res} (1933), \textsc{Zaremba} (1936),  \textsc{Pauc} (1936-41), \textsc{Ward} (1937),
\textsc{Saks} (1937), \textsc{Roger} (1938), \textsc{Choquet} (1943-48).}.

\section{Tangency}

The notion of tangency originated from geometric considerations in
antiquity. On the emergence of the coordinates of \textsc{Descartes}, analytic aspect
prevailed over the geometric view in tangency, also because of the growth of
infinitesimal calculus.

In \emph{Formulario Mathematico}\ \cite[(1908), p.\,313]{peano_1908}, a compendium
of mathematics known at the epoch, edited and mostly written by \textsc{Peano}\thinspace\footnote{%
In contrast to former versions that were written in French, the last (fifth)
version of \guillemotleft Formulario Mathematico\guillemotright\ (1908) was
written in \textquotedblleft latino sine flexione\textquotedblright.}, the
tangents of \textsc{Euclid} and \textsc{Descartes} are described in these terms:

\begin{quotation}
Euclide [...], dice que recta es tangente \guillemotleft $\varepsilon
\varphi \alpha \pi \tau \varepsilon \sigma \vartheta \alpha \iota $%
\guillemotright\ ad circulo [...] si habe uno solo puncto commune cum
circulo.

Nos pote applica idem Df [definition] ad ellipsi, etc.; sed non ad omni
curva.

Descartes, \emph{La G\'{e}om\'{e}trie }a.\,1637 \OE uvres, t.\,6, p.\,418
dice que tangente es recta que seca curva in duo puncto \guillemotleft ioins
en un\guillemotright ; id es, si \ae quatione que determina ce punctos de
intersectione habe duo \guillemotleft racines enti\`{e}rement \'{e}sgales%
\guillemotright .

Df [definition] considerato se transforma in P\textperiodcentered 0 [usual
definition], si nos considera per duo puncto \guillemotleft juncto in uno%
\guillemotright , ut limite de recto per duo puncto distincto.\thinspace\footnote{%
\translation{Euclid [...] says that a straight line is tangent to a circle [...] if it
has only one common point with the circle. One can apply the same definition
to an ellipse, and so on, but not to every curve. Descartes, \emph{La G\'{e}%
om\'{e}trie }a.\,1637 \OE uvres, t.\,6, p.\,418, says that a tangent is a
straight line that cuts a curve in two points \guillemotleft joined in one%
\guillemotright ; that is, the equation that determines these points of
intersection has two \guillemotleft entirely equal roots\guillemotright .
[The definition] considered [by Descartes] becomes [the usual definition] if
we mean by the points \guillemotleft joined in one\guillemotright\ the limit
of straight lines passing through two distinct points [when these tend to
one point].}}
\end{quotation}

A drawback of the predominance of analytic approach in geometry was that
tangency concepts were defined through an auxiliary system and not
intrinsically (that is, independently of a particular coordinate system).
Analytic approach to tangency requires that a figure, like a line or a
surface be defined via equations or parametrically, hence with the aid of
functions of some regularity. This constitutes another drawback, excluding,
for instance, figures defined by inequalities. On the other hand,
geometrically defined figures necessitate analytic translation before they
could be investigated for tangency.

The comeback to the geometric origin of tangency, and actually to a synergy
of both (geometric and analytic) aspects, is operated by the definitions of
tangency of arbitrary sets that use limits of homothetic figures. This
breakthrough was done by \textsc{Peano}  in \emph{Applicazioni
Geometriche} \cite[(1887)]{peano87}.

\emph{Synthetic geometry} started with \textsc{Euclid}, was axiomatized by \textsc{Pasch},
later by \textsc{Peano} and finally by \textsc{Hilbert}. \emph{Analytic geometry} (in the
original sense) was initiated by \textsc{Descartes} and enabled mathematicians to
reduce geometric problems to algebraic equalities, and thus to use algebraic
calculus to solve them. \emph{Vector geometry} of \textsc{Grassmann} potentiates the
virtues of both, synthetic and analytic, aspects of geometry.

In comparison with analytic methods, the classical geometric approach had
certainly an inconvenience of the lack of a system of standard operations
obeying simple algebraic rules, that is, of a calculus. In a letter of 1679
to \textsc{Huygens}, \textsc{Leibniz} postulated the need of a \emph{geometric calculus},
similarly to the already existing \emph{algebraic calculus}. This postulate
was realized by \textsc{Grassmann} in \emph{Geometrische Analyse} \cite[(1847)]%
{Grassmann:1847} and in \emph{Ausdehnungslehre} \cite[(1844, 1862)]{grass}. In \emph{%
Applicazioni Geometriche} \cite[(1887)]{peano87} \textsc{Peano} presented the geometric
calculus of \textsc{Grassmann} in order to treat geometric objects directly (without
coordinates), and in \emph{Calcolo Geometrico} \cite[(1888)]{calcolo88} refounded
the affine exterior algebra of \textsc{Grassmann} in three-dimensional spaces (see \textsc{Greco}, \textsc{Pagani}
\cite[(2009)]{greco_grass} for further details). In this way \textsc{Peano} eliminated the
inconvenience of the geometric approach mentioned above. This achievement
enabled him to develop a simple and sharp tangency theory abounding with
applications. Although \textsc{Peano}'s framework was that of 3-dimensional
Euclidean space, his method can be extended in an obvious way to arbitrary
dimensions (for example, the notion of angle between two subspaces can be
expressed in terms of the inner product multi-vectors).

\textsc{Peano}'s works permitted an easy access to the geometric calculus of 
\textsc{Grassmann} by the mathematical community at the end of 19th century\thinspace\footnote{%
See section 35: \emph{Begr\"{u}ndung der Punktrechnung durch G. Peano} in 
\cite[(1923)]{lotz} of the celebrated \emph{Encyklop\"{a}die der mathematischen
Wissenschaften}.}, in particular to the mathematicians of the Turin
University.

\section{Evolution of concepts of tangency in the work of Peano}

\label{evo} The interest of \textsc{Peano} in tangency goes back to 1882, two years
after he graduated from the university, when he discovered that the
definition of area of surface, given by \textsc{Serret} in his \emph{Cours de calcul
diff\'{e}rentiel et int\'{e}gral} \cite[p.\,293 (5th edition 1900)]{serret}
was defective. Indeed, \textsc{Serret} defined the area of a given surface as the
limit of the areas of polyhedral surfaces inscribed in that surface. \textsc{Peano}
found a sequence of polyhedral surfaces inscribed in a bounded cylinder so
that the corresponding areas tend to infinity \cite[(1902-1903), pp.\,300-301]%
{formulaire4}\thinspace \footnote{%
On reporting this discovery to his teacher \textsc{Genocchi}, \textsc{Peano} (24 years old) learned with
disappointment that \textsc{Genocchi} was already informed by \textsc{Schwarz} about the
defect of  \textsc{Serret}'s definition in 1882 (see \cite[p.\,9]{kennedy_corrected}%
).}. As \textsc{Peano} comments in that note

\begin{quotation}
On ne peut pas d\'{e}finir l'aire d'une surface courbe comme la limite de
l'aire d'une surface poly\'{e}drique inscrite, car les faces du poly\`{e}dre
n'ont pas n\'{e}cessairement pour limite les plans tangents \`{a} la surface.%
\thinspace\footnote{%
\translation{One cannot define the area of a curved surface as the limit of the area of
an inscribed polyhedral surface, because the faces of the polyhedron do not
necessarily tend to the tangent planes of that surface.}}
\end{quotation}

Lower (\ref{lower})  and upper (\ref{upper})  tangent cones constitute a final
achievement of \textsc{Peano}'s investigations started in \emph{Applicazioni
Geometriche }\cite[(1887)]{peano87}, where the lower tangent cone was already
defined explicitly as in (\ref{lower}), while the upper tangent cone was
implicitly used in \cite[(1887)]{peano87} in the proof of necessary optimality
conditions, and explicitly defined in \emph{Formulaire Math\'{e}matique} 
\cite[p.\,296]{formulaire4} of 1902-3 and in \emph{Formulario Mathematico} 
\cite[(1908) p.\,331]{peano_1908} as in (\ref{upper}). Apart from \cite[(1887)]{peano87} and \cite[(1908)]%
{peano_1908} \textsc{Peano} studies and uses tangency concepts in several other
works:  \emph{Teoremi su massimi e minimi
geometrici e su normali a curve e superficie}\}\cite[(1888)]{peano-mass}, \emph{%
Sopra alcune curve singolari} \cite[(1890)]{peano-sopra},  \emph{%
Elementi di calcolo geometrico} \cite[(1891)]{peano-elementi}, \emph{Lezioni di analisi infinitesimale }\cite[(1893)]{peano1893} and \emph{Saggio di
calcolo geometrico} \cite[(1895-96)]{saggio}.

Following this list we will trace the development of his ideas on tangency,
describing not only definitions and properties, but also his methods,
calculus rules and applications.

Peano managed to maintain exceptional coherence and precision during a
quarter of century of investigations on various and changing aspects of
tangency. Only a particular care, with which we perused his work, enabled us
to discern a couple of slight variations in the definitions, which, however,
did not induce \textsc{Peano} to any erroneous statement. For instance, \textsc{Peano} gives
an \emph{intrinsic definition} of tangent straight line to a curve, and also
another definition that is the tangent vector to the function representing
that curve. He underlines that the two notions are slightly different \cite[(1908),
p.\,332 (see properties P69.4, P70.1)]{peano_1908}

In \emph{Applicazioni Geometriche} \cite[(1887)]{peano87}, after having presented
elements of the geometric calculus of \textsc{Grassmann} (point, vector, bi-vector,
tri-vector\thinspace\footnote{%
A bi-vector is the exterior product of 2 vectors, a tri-vector is the
exterior product of 3 vectors. Vectors, bi-vectors and tri-vectors are used
by \textsc{Peano} in 1888 in replacement of the corresponding terms of segment, area
and volume adopted in \emph{Applicazioni Geometriche} \cite[(1887)]{peano87}.},
scalar product and linear operations on them), \textsc{Peano} defines limits of
points and vector-type objects (vectors, bi-vectors, tri-vectors) and proves
the continuity and differentiability of the operations of addition, scalar multiplication, scalar
 product and products of vectors (see pages 39--56 of \emph{Applicazioni Geometriche}). 

Moreover he defines limits of straight lines and of planes. 
Straight lines and planes are seen
by \textsc{Peano} as sets of points, so that their limits are instances of a general
concept of convergence of variable sets: the \emph{lower limit} \eqref{lower}.  Accordingly, a variable straight line (a
variable plane) $A_{t}$ converges to a straight line (plane) $A$ as a parameter $t$ tends to some finite or infinite quantity, if
\begin{equation}
A\subset \limfunc{Li}\nolimits_{t}A_{t},  \label{Liminf}
\end{equation}%
 that is, if the distance $d(x,A_{t})$ converges to $0$ for each $x\in A$.
Then he checks meticulously (without using coordinates) the continuity of various relations involving
points, straight lines and planes. For instance,

\begin{itemize}
\item[$(i)$] A variable straight line $L_{t}$ converges to a straight line $%
L $ if and only if for two distinct points $x,y\in L$ the distances of $%
L_{t} $ to $x $ and $y$ tend to $0$.

\item[$(ii)$] A variable plane $P_{t}$ converges to a plane $P$ if and only
if for non-colinear points $x,y,z\in P$ the distances of $P_{t}$ to $x,y$
and $z$ tend to $0$.

\item[$(iii)$] If two variable straight lines $L_{t}$ and $M_{t}$ converge to the non-parallel straight lines $L$ and $M$, respectively, then the straight line $N_{t}$ which meets perpendicularly both  $L_{t}$ and $M_{t}$, converge to   the straight line $N$ which meets perpendicularly both  $L$ and $M$. 
\end{itemize}

In \emph{Applicazioni Geometriche }\cite[(1887), p.\,58]{peano87} \textsc{Peano}
defines 

\begin{definition}
\label{def:straight} A \emph{tangent straight line} of a \emph{curve} $C$ at
a point $x\in C$ is the limit of the straight line passing through $x$ and
another point $y\in C$ as $y$ tends to $x$.
\end{definition}

For \textsc{Peano}, a \emph{curve} $C$ is a subset of Euclidean space such
that $C$ is homeomorphic to an interval $I$ of the real line, so that $%
C=\left\{ C(t):t\in I\right\} $ can be seen as depending on a parameter $%
t\in I$.  He gives a description of the tangent straight line in the case where the
derivatives $C^{(k)}(\hat{t})$ are null for $k<p$ and $C^{(p)}(\hat{t})\neq
0 $. Moreover,

\begin{proposition}
\label{prop:Ppara}\cite[(1887), teorema II, p.\,59]{peano87} If $C$ is continuously
differentiable and $C^{\prime }(\hat{t})\neq 0$, then the tangent straight
line $L$ is the limit of the lines passing through $x,y\in C$ as $x,y$ tend to $C(\hat{t})$ and $x\neq
y$.
\end{proposition}

Notice that Proposition \ref{prop:Ppara}  makes transparent the relation between
paratangency and the continuity of derivative (see \eqref{paradiff} for a sequential description of paratangent vector). Paratangency to curves and surfaces was used by \textsc{Peano} also in other instances in \emph{Applicazioni geometriche} \cite[(1887), p.\,163, 181-184]{peano87} to evaluate the infinitesimal quotient of the length of an arc and its segment or its projection.

After a study of mutual positions of a curve and its tangent straight lines,
\textsc{Peano} gives rules for calculating the tangent straight line to the graph of a
function of one variable and to a \emph{curve} given by an equation $%
f(x,y)=0 $, or by two equations%
\begin{equation*}
f(x,y,z)=0\text{ and }g(x,y,z)=0,
\end{equation*}%
for which he needs the implicit function theorem. Incidentally, he
presented, for the first time in 1884 in a book form \cite{Genocchi}, the
implicit function theorem proved by \textsc{Dini} in 1877-78 in his lectures \cite%
[pp.\,153-207]{Dini} and provided a new proof, much shorter than the original
demonstration of \textsc{Dini}.

\textsc{Peano} gave numerous examples of application of these calculus rules, among
others, to parabolas of arbitrary order, logarithmic curve, Archimedean
spiral, logarithmic spiral, concoids (e.g., limacon of Pascal, cardioid),
cissoids (e.g., lemniscate).

Successively \textsc{Peano} defines

\begin{definition}
\label{def:plane} A \emph{tangent plane} to a surface $S$ at a given point $%
x\in S$ is the plane $\alpha$ such that the acute angle between $\alpha$ and each
straight line passing through $x$ and another point $y\in S$ tends to $0$ as 
$y$ tends to $x$.
\end{definition}

A \emph{surface} is assumed to be a subset (of Euclidean space)
homeomorphic to a rectangle. Several properties of tangent planes are then
proved intrinsically, by geometric calculus, without the use of coordinates
or parametric representations.

He also calculates intrinsically the tangent planes of many classical
surfaces, like cones, cylinders and revolution figures, and more generally,
surfaces obtained by a rigid movement of a curve. As he did before with
curves, \textsc{Peano} calculates tangent planes to the graphs of functions of two
variables as well as to surfaces given by equations and parametrizations. As
for curves, he gives analytic criteria on the position of a surface with
respect to its tangent planes.

The novelty does not consist of a description of particular cases of
tangency, but of the precision and the refinement of the analysis of
conditions that are necessary for tangency, which characterize the methods
of geometric calculus.

\section{Remarks on relationship between tangency and differentiability}

Most sophisticated examples of calculation of tangent planes come from
geometric operations, like geometric loci (described in terms of distance
functions from points, straight lines and planes). They are based on the
notion of differentiability introduced by \textsc{Peano} (called nowadays \emph{Fr%
\'{e}chet differentiability}). An essential tool is the following theorem on
differentiability of distance functions
\footnote{A detailed study of regularity of distance function was carried out for the first time  by \textsc{Federer}
in \cite[(1959)]{Federer}.
}.

\begin{theorem}
\cite[(1887), pp.\,139-140]{peano87} Let $F$ be a subset of the Euclidean space $X$
such that there exists a continuous function $\gamma :X\rightarrow F$ so
that $d(x,\gamma (x))=d(x,F)$. Then the distance function $x\mapsto d(x,F)$
is differentiable at each point $\hat{x}\notin F$ and the derivative is
equal to $\dfrac{\hat{x}-\gamma (\hat{x})}{\Vert \hat{x}-\gamma (\hat{x}%
)\Vert }$.
\end{theorem}

Finally, in the last chapter of \emph{Applicazioni Geometriche}, \textsc{Peano}
introduces lower   affine tangent cone of an arbitrary  subset of the Euclidean space $X$ \cite[(1887), p.\,305]{peano87}. The \emph{%
lower affine tangent cone }$\func{tang}(F,x)$ of $F$ at $x$ (for arbitrary $%
x\in X$) is given by the blowup 
\begin{equation}
\func{tang}(F,x)=\limfunc{Li}\nolimits_{h\rightarrow +\infty }\left(
x+h(F-x)\right),\thinspace\footnote{%
Observe that the lower affine tangent cone is an affine version of the lower
tangent cone, since $\func{tang}(F,x)=x+\limfunc{Tan}\nolimits^{-}(F,x)$.}
\label{tang1}
\end{equation}%
hence, by \eqref{lowerlimit} 
\begin{equation}
y\in \func{tang}(F,x)\Longleftrightarrow \lim\nolimits_{t\rightarrow 0^{+}}%
\tfrac{1}{t}d(x+t(y-x),F)=0.  \label{tang}
\end{equation}%
\textsc{Peano} claims that $\func{tang}(F,x)$ \textquotedblleft
generalizes\textquotedblright\ the tangent straight line of a curve and the
tangent plane of a surface. Actually, there is a discrepancy between (\ref%
{tang}) and Definitions \ref{def:straight} and \ref{def:plane}, because the
tangent defined above is a cone that need not be a straight line (resp. a
plane).\thinspace \footnote{%
If $F:=\left\{ (x,y)\in \mathbb{R}^{2}:y=\sqrt{\left\vert x\right\vert }%
\right\} ,$ then the tangent straight line to $F$ at the origin in the sense
of Definition \ref{def:straight} is $\left\{ (x,y)\in \mathbb{R}%
^{2}:x=0\right\} $, while $\func{tang}(F,(0,0))=\left\{ (x,y)\in \mathbb{R}%
^{2}:x=0,y\geq 0\right\} $.}

Tangency was principally used by \textsc{Peano} for the search of maxima and
minima with the aid of necessary conditions of optimality. Many of
optimization problems considered in \emph{Applicazioni Geometriche} are
inspired by geometry, for example: \textquotedblleft Find a point that
minimizes the sum of the distances from given three
points\textquotedblright\ \cite[(1887), p.\,148]{peano87}.

Necessary optimality conditions (see Theorem \ref{thm:regula} below) given
in \emph{Applicazioni Geometriche}, reappear in \emph{Formulario Mathematico}
formulated with the aid of the upper affine tangent cone. The \emph{upper
affine tangent cone} is defined by the blowup 
\begin{equation}
\func{Tang}(F,x)=\limfunc{Ls}\nolimits_{h\rightarrow +\infty }\left(
x+h(F-x)\right) .\thinspace\footnote{%
Observe that the upper affine tangent cone is an affine version of the upper
tangent cone, since $\func{Tang}(F,x)=x+\limfunc{Tan}\nolimits^{+}(F,x)$.}
\label{tang2}
\end{equation}

Hence, by (\ref{upperlimit}),\begin{equation}
y\in \limfunc{Tang}(F,x)\Longleftrightarrow \lim \inf\nolimits_{t\rightarrow
0^{+}}\tfrac{1}{t}d(x+t(y-x),F)=0.
\end{equation}

\begin{theorem}[\textsc{Peano}'s Regula]
\label{thm:regula}If $f\colon \mathbb{R}^{n}\rightarrow \mathbb{R}$ is
differentiable at $x\in A\subset \mathbb{R}^{n}$ and $f(x)=\max \left\{
f(y):y\in A\right\} $, then%
\begin{equation}
\langle Df(x),y-x\rangle \leq 0\text{ for each }y\in \limfunc{Tang}(A,x),
\label{max=cond}
\end{equation}%
where $Df(x)$ denotes the gradient of $f$ at $x$.
\end{theorem}

This theorem was formulated in \emph{Formulario Mathematico} \cite[(1908), p.\,335]%
{peano_1908} exactly as above, but was proved informally already in \emph{%
Applicazioni Geometriche} \cite[(1908), p.\,143-144]{peano87} (without an explicit
definition of the upper affine tangent cone). Condition (\ref{max=cond}) is
best possible in the following sense:\thinspace \footnote{%
Indeed, if $w\in \limfunc{Nor}(A,x)$ then we define $f:\mathbb{R}%
^{n}\rightarrow \mathbb{R}$ as follows $f(y)=\langle w,y\rangle $ for each $%
y $ with the exception of $y\in A\cap \left\{ y:\langle w,y-x\rangle \geq
0\right\} $, for which $f(y)=f(x)$.} 
\begin{equation}
\left\{ Df(x):f\text{ is differentiable at }x\text{ and }\max%
\nolimits_{A}f=f(x)\right\} =\limfunc{Nor}(A,x),  \label{best}
\end{equation}%
where the usual \emph{normal cone }(defined by \textsc{Federer} \cite{Federer} in
1959) is 
\begin{equation*}
\limfunc{Nor}(A,x):=\left\{ w\in \mathbb{R}^{n}:\langle w,y-x\rangle \leq 0%
\text{ for each }y\in \limfunc{Tang}(A,x)\right\} .
\end{equation*}

The equivalence of differentiability and of the existence of tangent
straight line was considered as evident from the very beginning of
infinitesimal calculus.

In case of functions of several variables however relationship between
differentiability and tangency remained vague, partly because the very
notion of tangency was imprecise. 

Ways to a definition of tangency were
disseminated with pitfalls as witness several unsuccessful attempts.
For instance, \textsc{Cauchy} confused partial differentiability and differentiability, that is,
the existence of total differential\thinspace\footnote{%
Also the relation between separate and joint continuity was elucidated
long after erroneous claims of \textsc{Cauchy} in 1821 in \cite{cauchy}. A classical example of function of two variables that is separately
continuous but not continuous was provided by \textsc{Peano} in \cite[(11884) p.\,173]%
{Genocchi}: $(x,y)\mapsto{xy}/{(x^{2}+y^{2})}. $}. \textsc{Thomae} was the first to distinguish the
two concepts in \cite[(1875), p.\,36]{thomae} by supplying simple counter-examples.

Differentiability of a function of several variables was defined by \textsc{Peano} in 
\cite[(1887)]{peano87}, as it is defined today under the name of \emph{Fr\'{e}chet
differentiability} and reappears in his \emph{Formulario Mathematico} in 
\cite[(1908) p.\,330]{peano_1908}. With the exception of \cite[(1891), p.\,39]{peano-elementi},
where he observes that the existence of total differential could be taken as
a definition of differentiability, \textsc{Peano} uses, in numerous applications, the
continuity of partial derivatives, which amounts to strict
differentiability. He notices in \cite[(1892)]{peano_mathesis} that strict
differentiability is equivalent to the uniform convergence of the difference
quotient to the derivative, as he also does in an epistolary exchange (see 
\cite[(1884)]{peano-lettre} and \cite[(1884)]{peano-reponse}), concerning the hypotheses of
the mean value theorem in the book of \textsc{Jordan} \cite[(1882)]{Jordan}\thinspace\footnote{%
\textsc{Peano} points out that it is enough to assume differentiability, and not
continuous differentiability as did \textsc{Jordan} and \textsc{Cauchy}.}. The idea of strict
differentiability is extended by \textsc{Peano} in a spectacular way to the theory of
differentiation of measures (see \textsc{Greco}, \textsc{Mazzucchi} and \textsc{Pagani} \cite{greco-sonia} for details).

\textsc{Peano} criticizes various existent definitions of tangency  \cite[(1908) p.\,333]{peano_1908}:

\begin{quotation}
Plure Auctore sume ce proprietate ut definitione. \guillemotleft Plano
tangente ad superficie in suo puncto $p$\guillemotright\ es definito ut 
\guillemotleft plano que contine recta tangente in $p$ ad omni curva,
descripto in superficie, et que i trans $p$\guillemotright.\thinspace\footnote{%
\translation{Several authors take this property as a definition: \guillemotleft a
tangent plane to a surface at its point $p$\guillemotright\ is defined as 
\guillemotleft a plane that contains the tangent straight line at $p$ of
every curve traced on the surface and passing through $p$\guillemotright.}}
\end{quotation}

As counter-examples to this definition, \textsc{Peano} quotes a logarithmic spiral at
its pole\thinspace\footnote{%
called also a \emph{miraculous spiral} (\emph{spira mirabile} in latino sine
flexione), after the Latin name \emph{spira mirabilis }given to it by \textsc{J.
Bernoulli}, that is, a curve described in polar coordinates $(r,\theta )$ by $%
r=ae^{b\theta }$. 
The pole is the
origin of $\mathbb{R}^{2}$.} and a loxodrome at its poles.
He continues

\begin{quotation}
Aliquo Auctore corrige pr\ae cedente, et voca plano tangente \guillemotleft %
plano que contine tangente ad dicto curvas, que habe tangente\guillemotright %
.\thinspace\footnote{%
\translation{Other authors correct the preceding [definition], and call a tangent plane 
\guillemotleft the plane that contains the tangent to those [said] curves
that have a tangent [straight line]\guillemotright .}}
\end{quotation}

He constructs a counter-example\thinspace\footnote{%
By rotating around the $x$-axis in the space of $(x,y,z)$, the function%
\begin{equation*}
y=\left\{ 
\begin{array}{c}
x\sin \left( \frac{1}{x}\right) \text{ if }x\neq 0 \\ 
0\text{ if }x=0%
\end{array}%
\right. ,
\end{equation*}%
that he had introduced. Recall that at that epoch, a \emph{curve} is assumed
to be continuous.} to this definition that was adopted, among others, by
\textsc{Serret} \cite[p.\,370]{serret}. \textsc{Bertrand}, one of the most famous and
influential French mathematicians of 19th century, writes in \cite[(1865), p.\,15]%
{bertrand}

\begin{quotation}
Le plan tangent d'une surface en un point est le plan qui, en ce point,
contient les tangentes \`{a} toutes les courbes trac\'{e}es sur la surface.
\thinspace\footnote{%
\translation{The tangent plane of a surface at a point is the plane that, at this point,
includes all the tangent lines to all the curves drawn on the surface.}}
\end{quotation}

The literature abounds with observations, mostly in view of didactic use, on
the relation between the notion of tangent plane at the graph of the function%
\begin{equation}
z=f(x,y)  \label{sureq}
\end{equation}%
and the differentiability of $f$ at interior points of the domain of $f$.
For example, in \cite[(1911)]{Frechet_diff} \textsc{Fr\'{e}chet} observes\thinspace \footnote{%
\translation{A function $f(x,y)$ has a differential in my sense at $(x_{0},y_{0})$, if
the surface $z=f(x,y)$ admits at this point a unique tangent plane
non-parallel to $Oz$: $z-z_{0}=p(x-x_{0})+q(y-y_{0}$). Then this
differential is, by definition, 
\begin{equation}
p\Delta x+q\Delta y,  \tag{*}
\end{equation}%
where $\Delta x$, $\Delta y$ are arbitrary increments of $x$, $y$. [\dots ]
The analytic form of this definition is the following: [\dots ] A function $%
f(x,y)$ has a differential in my sense at $(x_{0},y_{0})$ if there exists a
linear homogeneous function (*) of increments that differs from $\Delta f$
[\dots ] by an infinitesimal with respect to the distance $\Delta $ of the
points $(x_{0},y_{0})$, $(x_{0}+\Delta x,y_{0}+\Delta y).$}}

\begin{quotation}
Une fonction $f(x,y)$ a une diff\'{e}rentielle \`{a} mon sens au point $%
(x_{0},y_{0})$, si la surface $z=f(x,y)$ admet en ce point un plan tangent
unique non parall\`{e}le \`{a} $Oz$: $z-z_{0}=p(x-x_{0})+q(y-y_{0}$). Et
alors cette diff\'{e}rentielle est par d\'{e}finition l'expression 
\begin{equation}
p\Delta x+q\Delta y,  \tag{*}
\end{equation}%
o\`{u} $\Delta x$, $\Delta y$ sont des accroissements arbitraires de $x$, $y$%
. [\dots ] La forme analytique de cette d\'{e}finition est la suivante:
[\dots ] Une fontion $f(x,y)$ admet une diff\'{e}rentielle \`{a} mon sens au
point $(x_{0},y_{0})$ s'il existe une fonction lin\'{e}aire et homog\`{e}ne
(*) des accroissements, qui ne diff\`{e}re de l'accroissement $\Delta f$
[\dots ] que d'un infiniment petit par rapport \'{a} l'\'{e}cart $\Delta $
des points $(x_{0},y_{0})$, $(x_{0}+\Delta x,y_{0}+\Delta y)$,\thinspace\footnote{%
Fr\'{e}chet forgets that in order that a tangent plane imply
differentiability, it is necessary to assume the continuity of $f$ at $%
(x_{0},y_{0})$.}
\end{quotation}

Surely, this definition would be certainly more precise if \textsc{Fr\'{e}chet} had
defined his concept of tangency
\footnote{In \cite[(1964), p.\,189]{frechet64} \textsc{Fr\'{e}chet} gives the following definition of the tangent plane that slightly differs from that of \textsc{Bertrand}: 
\begin{quotation}
Pr{\'e}cisons d'abord que nous entendons par plan tangent {\`a} [une surface] $S$ au point $(a,b,c)$ un plan qui soit lieu des tangentes aux courbes situ{\'e}es sur $S$ et passant par ce point (s'entendant de celles de ces courbes qui ont effectivement une tangente en ce point).
\end{quotation}
\translation{Let us first make precise that by tangent plane to [a surface] $S$ at a point $(a,b,c)$, we mean a plane that is the locus of tangent lines to the curves lying on $S$ and passing through this point
(that is, to those curves that have effectively a tangent line at that point).} }.

\textsc{Wilkosz} characterizes in \cite[(1921)]{wilkosz} differentiability in terms of
non-vertical tangent half-lines that form a single plane and are uniform
limits of the corresponding secants. It is notable that he acknowledges
\textsc{Stolz} and \textsc{Peano} as creators of the notion of total differential.

\textsc{Saks} defines in \cite[(1933)]{saks33} differentiability as the existence of a
tangent plane at (\ref{sureq}) in the sense of Definition \ref{def:plane}.
Consequently, a tangent plane of \textsc{Saks} can contain  vertical lines.

\textsc{Tonelli} defines in \cite[(1940)]{tonelli} differentiability as the existence of a
tangent plane in the sense of Definition \ref{def:plane}, provided that the
orthogonal projection of (\ref{sureq}) on the tangent plane is open at the
point of tangency. His notion of differentiability coincides with
the modern concept of differentiability.
\section{Characterizations of differentiability} \label{diff-cara}

\textsc{Guareschi} and \textsc{Severi} characterized differentiability in terms of tangency of
their graphs (for functions defined on subsets of Euclidean space). At the
same period also \textsc{Bouligand} studied tangency, but his perception of the
relationship between differentiability and tangent cones remained vague \cite%
[(1932), pp.\,68-71]{bouligand}.

\textsc{Guareschi} and \textsc{Severi} stress that the originality of their approach consists
in defining a \emph{total differential} of a function $f$ defined on an 
\emph{arbitrary} subset $A$ of Euclidean space at an accumulation point of 
$A$. Consequently, their definition cannot hinge on traditional \emph{%
partial derivatives}. In \cite[(1934)]{Guareschi-diff}, \textsc{Guareschi}, using a notion of 
\emph{tangent figure} of \textsc{Severi} \cite[(1929, 1931)]{Severi-para,Severi-Krak},
introduces a \emph{linear tangent space} in order to characterize existence
and uniqueness of total differentials. Both refer to the notion of
differentiability of \textsc{Stolz} \cite[(1893)]{stolz}.

The \emph{tangent figure} of \textsc{Severi} is defined (only at accumulation points)
as the union of all tangent half-lines (that he called \emph{semi-tangents}%
), in the same way as \textsc{Saks} describes in \cite[(1933), p.\,262]{Saks} the \textsc{Bouligand} 
\emph{contingent cone} \cite{bouligand}. As observed in \cite[p.\,501]%
{DG-Peano_top}, \textsc{Severi}'s tangent figure is precisely the \emph{upper tangent
cone} (\ref{upper}) of \textsc{Peano}; as we have already noted, although \textsc{Severi}
cites \textsc{Bouligand} and \textsc{Saks}, he never quotes \textsc{Peano} (see footnote \ref{ftnote}).
Nevertheless in \cite[p.\,23 (footnote)]{Severi-Poinc} \textsc{Severi} writes in 1949

\begin{quotation}
[...] nostro grande logico matematico Giuseppe Peano, che fu mio maestro ed
amico e della cui intuizione conobbi tutta la forza.\thinspace\footnote{%
\translation{[\dots] our great logician and mathematician Giuseppe Peano, who was my mentor
and friend, of whose intuition I knew all the strength.}}
\end{quotation}

As we mentioned, neither \textsc{Guareschi} cited \textsc{Peano}. He however did not forget to
send the following telegram on the 70th birthday of \textsc{Peano}.

\begin{quotation}
Esprimo illustre scienziato ammirazione augurio lunga feconda attivit\`{a}.%
\thinspace\footnote{%
\translation{I express, illustrious scientist, admiration [and] wishes of long [and]
fertile activity.}}
\end{quotation}

\textsc{Guareschi} \cite[(1934), p.\,177]{Guareschi-diff} reformulates the \textsc{Severi}'s definition
of \emph{upper affine tangent cone} with the aid of conical neighborhoods.
If $\hat x$ is a point and $h$ is a non-zero vector of Euclidean space, then a 
\emph{conical neighborhood} $C(\hat x,h,r,\alpha )$ of a half-line, starting at $%
\hat x $ in the direction $h$, is the intersection of a sphere (of a radius $r>0$%
) centered at $\hat  x$ with a revolution cone of solid angle $\alpha $ around the
axis $h$. A half-line at $\hat x$ in the direction $h$ is \emph{tangent} to $A$
at $\hat x$ if and only if $C(\hat x,h,r,\alpha )\cap A\setminus \left\{ \hat x\right\}
\neq \varnothing $ for every $r>0$ and $\alpha >0$.

In fact, this definition had been already given by \textsc{Cassina} in \cite[(1930)]%
{Cassina-Milano}. \textsc{Cassina} presented it as an alternative description of the%
\emph{\ lower tangent cone} (\ref{lower}) from \emph{Applicazioni Geometriche%
}; \textsc{Cassina}'s definition is however equivalent to the \emph{upper tangent cone%
} (\ref{upper}), for which \textsc{Cassina} proves  the following new fact%
\thinspace\footnote{%
We regret to have forgot to cite in \cite{DG-Peano} this contribution of
\textsc{Cassina}, which is parallel to those of \textsc{Bouligand} and \textsc{Severi}.} that includes
a later result of \textsc{Severi} \cite[(1931)]{Severi-Krak}.

\begin{theorem}[{\textsc{Cassina} \cite[(1930)]{Cassina-Milano}}]
\label{thm:Cassina} There exists a tangent half-line of $A$ at $\hat x$ if and
only if $\hat x$ is an accumulation point of $A$.
\end{theorem}

\textsc{Guareschi}'s characterization of differentiability is as follows. By $%
\limfunc{graph}(f)$ we denote the \emph{graph} of a function $f:A\rightarrow 
\mathbb{R}$, where $A\subset \mathbb{R}^{n}$. Of course, a hyperplane $H$ in 
$\mathbb{R}^{n}\times \mathbb{R}$ is a graph of an affine function from $\mathbb{R}^{n}$ to $\mathbb{R}$, whenever $%
H$ does not include  vertical lines.
 
\begin{theorem}[{\textsc{Guareschi} \cite[(1934), p.\,181]{Guareschi-diff}}]
\label{thm:Guareschi} Let $A\subset \mathbb{R}^{n}$ and let $\hat x\in A$ be an
accumulation point  of $A$. A function $f:A\rightarrow \mathbb{R}$,
continuous at $\hat x$, is differentiable at $\hat x$ if and only if $\limfunc{Tan}%
^{+}(\limfunc{graph}(f),(\hat x,f(\hat x)))$ is included in a hyperplane without 
vertical lines.
\end{theorem}

The \emph{linear tangent space} of \textsc{Guareschi} at an accumulation point $\hat x$ of $A$ is exactly the affine
space spanned by the upper affine tangent cone of $A$ at $\hat x$; its dimension is called by \textsc{Guareschi}, \emph{accumulation dimension}\emph{\ } of $A$ at point $\hat x$  \cite[(1934), p.\,184]{Guareschi-diff}.

The \emph{total differential} of a function $f:A\rightarrow \mathbb{R}$ at
an accumulation point $\hat x$ of $A$ with $\hat x\in A$ is defined as a linear map $L:\mathbb{R}%
^{n}\rightarrow \mathbb{R}$ such that%
\begin{equation*}
\lim\nolimits_{A\ni y\rightarrow x}\dfrac{\left\vert
f(y)-f(\hat x)-L(y-\hat x)\right\vert }{\left\Vert y-\hat x\right\Vert }=0.
\end{equation*}%
Using these notions, \textsc{Guareschi} reformulates Theorem \ref{thm:Guareschi}:

\begin{theorem}[{\textsc{Guareschi} \cite[(1934), p.\,183]{Guareschi-diff}}]
\label{thm:guareschi2} Let $f$ be a real function on a subset of Euclidean
space of dimension $n$. If $\limfunc{Tan}^{+}(\limfunc{graph}(f),(\hat x,f(\hat x)))$
does not include  vertical lines, then the following  properties hold:

\begin{enumerate}
\item there exists a total differential of $f$ at $\hat x$ if and only if the accumulation
dimension of $\limfunc{graph}(f)$ at $(\hat x,f(\hat x))$ is not greater than $n$;

\item a total differential of $f$ at $\hat x$ is unique if and only if the accumulation
dimension of $\limfunc{graph}(f)$ at $(\hat x,f(\hat x))$ is $n$.
\end{enumerate}
\end{theorem}

Therefore there is a one to one correspondence between total differentials
and hyperplanes without vertical lines that include the tangent figure $%
\limfunc{Tan}^{+}(\limfunc{graph}(f),(\hat x,f(\hat x)))$.

\textsc{Severi} presented the paper \cite[(1934)]{Guareschi-diff} of \textsc{Guareschi} to the \emph{Reale
Accademia d'Italia} on the 10th November 1933, having suggested to the author
several simplifications and generalizations. Subsequently, \textsc{Severi}
reconsidered the topic in \cite[(1934)]{Severi-diff} and extended the results of 
\textsc{Guareschi}; he presented in a clear way the ideas of \textsc{Guareschi},
which originally were introduced with complex technicalities. 

The
differentiability results of \cite[(1934)]{Severi-diff} can be restated (and
partially reinforced) in the following, more modern way.

\begin{theorem}[\textsc{Severi-Guareschi}]
\label{thm:sintesi} Let $f:A\rightarrow \mathbb{R}^{k}$ where $A\subset 
\mathbb{R}^{m}$, and let $\hat{x}\in A$ be an accumulation point of $A$. Let $L:$ 
$\mathbb{R}^{m}\rightarrow \mathbb{R}^{k}$ be a linear map. Then the
following properties are equivalent:

\begin{enumerate}
\item $f$ is differentiable at $\hat{x}$ and $L$ is a total differential of $%
f$ at $\hat{x}$;

\item $f$ is continuous at $\hat{x}$ and $\limfunc{Tan}\nolimits^{+}\left( 
\limfunc{graph}(f),(\hat{x},f(\hat{x}))\right) \subset\limfunc{graph}(
L) $;

\item $\lim\nolimits_{n}\dfrac{f(x_{n})-f(\hat{x})}{\Vert x_{n}-\hat{x}\Vert 
}=L(v)$ for each $v\in \mathbb{R}^{m}$ and for every sequences $\left\{
x_{n}\right\} _{n}\subset A$ such that $\lim\nolimits_{n}x_{n}=\hat{x}$ and
$\lim\nolimits_{n}\dfrac{x_{n}-\hat{x}}{\Vert x_{n}-\hat{x}\Vert }=v$;

\item $L(v)=\lim\nolimits_{\substack{ w\rightarrow v  \\ t\rightarrow 0^{+}}}%
\dfrac{f(\hat{x}+tw)-f(\hat{x})}{t}$ for every $v\in \limfunc{Tan}%
\nolimits^{+}(A,\hat{x})$.
\end{enumerate}
\end{theorem}

Condition (2) of the theorem above encompasses Theorem \ref{thm:guareschi2}.
Condition (3) corresponds to \cite[(1934), pp.\,183-184]{Severi-diff} of \textsc{Severi}.
Condition (4) represents the total differential in terms of the directional
derivatives along tangent vectors \cite[(1934), p.\,186]{Severi-diff}, called \emph{%
perfect derivatives} by \textsc{Guareschi} \cite[(1934) p.\,201]{Guareschi-diff}. These
derivatives are usually formulated in terms of (just mentioned) conical
neighborhoods, and called \emph{Hadamard derivatives}.\thinspace\footnote{%
In spite of our efforts, we were unable to find these derivatives in
\textsc{Hadamard}'s papers. The reference \cite[(1923)]{hadamard} usually mentioned in this context
does not contain any pertinent  fact.}

Another condition equivalent to those of Theorem \ref{thm:sintesi} turns out
to be very instrumental in effective calculus of total differential%
\thinspace\footnote{%
Because of his pedagogical experience, in which the condition was frequently
of great help, the second author named it the \emph{Cyrenian Lemma},
referring to \textsc{Simon of Cyrene} who helped to carry the \textsc{Christ}'s cross.}.

\begin{proposition}
(Cyrenian Lemma) A function $f$ is differentiable at $\hat{x}$ and $L$ is a
total differential of $f$ at $\hat{x}$ if and only if $\lim\nolimits_{n}%
\dfrac{f(x_{n})-f(\hat{x})}{\lambda _{n}}=L(v)$ for each $v\in \mathbb{R}%
^{m} $ and for every sequences $\left\{ x_{n}\right\} _{n}\subset A$ and $%
\left\{ \lambda _{n}\right\} _{n}\subset \mathbb{R}_{++}$ such that $%
\lim\nolimits_{n}\lambda _{n}=0,\lim\nolimits_{n}x_{n}=\hat{x}$ and
$\lim\nolimits_{n}\dfrac{x_{n}-\hat{x}}{\lambda _{n}}=v$.\thinspace\footnote{%
As an instance of its usefulness, let us calculate the total differential at 
$(0,0)$ of%
\begin{equation*}
f(x,y):=x+y+\sqrt[2]{y^{3}(x-y)^{3}},\;\limfunc{dom}f:=\left\{ (x,y)\in 
\mathbb{R}^{2}:y^{3}(x-y)^{3}\geq 0\right\} ,
\end{equation*}%
that was calculated (over several pages) by \textsc{Guareschi} in \cite[(1934), p.\,190-194]%
{Guareschi-diff}. In fact if $\lambda _{n}\rightarrow 0^{+},\limfunc{dom}%
f\ni (x_{n},y_{n})\rightarrow (0,0)$ and $\dfrac{1}{\lambda _{n}}\left[
(x_{n},y_{n})-(0,0)\right] \rightarrow (v,w)$, then the function $L:\mathbb{R}^{2}\to\mathbb{R}$ is well defined by
\begin{equation*}
L(v,w):=\lim\nolimits_{n}\dfrac{1}{\lambda _{n}}\left[ f(x_{n},y_{n})-f(0,0)%
\right] =\lim\nolimits_{n}\left( \dfrac{x_{n}}{\lambda _{n}}+\dfrac{y_{n}}{%
\lambda _{n}}+\sqrt[2]{\dfrac{y_{n}^{2}}{\lambda _{n}^{2}}%
y_{n}(x_{n}-y_{n})^{3}}\right) =v+w.
\end{equation*}%
ant it is linear. Hence, by \emph{Cyrenian Lemma}, $L$ is a total differential of $f$ at $(0,0)$.}
\end{proposition}

Theorem \ref{thm:sintesi} reformulates certain ingredients of the
characterizations above in a (hopefully) comprehensive way. For instance,
the non-verticality condition is incorporated in each of the conditions
(2-4). It is worthwhile to make explicit the particular case of
differentiability at interior points of the domain.

\begin{proposition}
\label{prop:diff=int} Let $A\subset \mathbb{R}^{m}$ and let $\hat{x}\in 
\limfunc{int}A$. A map $f:A\rightarrow \mathbb{R}^{k}$ is differentiable at $%
\hat{x}$ if and only if

\begin{enumerate}
\item $f$ is continuous at $\hat{x}$;\label{C1}

\item For each $v\in \mathbb{R}^{m}$ the directional derivative $\dfrac{%
\partial f}{\partial v}(\hat{x})$ exists and is linear in $v$;\label{C2}

\item $\limfunc{Tan}\nolimits^{+}(f,(\hat{x},f(\hat{x}))$ is a vector space
of dimension $m$.\label{C3}
\end{enumerate}
\end{proposition}

Observe that Condition (2) is usually referred to as \emph{G\^{a}teaux
differentiability}. In Proposition \ref{prop:diff=int} above none of the
three conditions can be dropped.

\begin{example}
Let $m:=2,k:=1,A:=\mathbb{R}^{2},\hat{x}:=(0,0).$

\begin{enumerate}
\item $f(x,y):=\left\{ 
\begin{array}{c}
1\text{ if }y=x^{2}\neq 0 \\ 
0\text{ otherwise}%
\end{array}%
\right. $ fulfills (\ref{C2}) and (\ref{C3}) but does not fulfill (\ref{C1}).

\item $f(x,y):=\sqrt[3]{x}$ fulfills (\ref{C1}) and (\ref{C3}) but not (\ref%
{C2}).

\item $f(x,y):=\left\{ 
\begin{array}{c}
x\text{ if }y=x^{2} \\ 
0\text{ otherwise}%
\end{array}%
\right. $ fulfills (\ref{C1}), (\ref{C2}) but not (\ref{C3}).
\end{enumerate}
\end{example}

\section{Characterizations of strict differentiability} \label{strictdiff-cara}

Till the installation of the today concept of differentiability, the
continuity of partial derivatives had been used to affirm the existence of
total differential. As it turned out that this condition is sufficient but
not necessary, \textsc{Severi} wanted to find an additional property of the total
differential corresponding to the continuity of partial derivatives. He
discovered that, for the internal points of the domain, \emph{strict
differentiability} (\ref{strict differ}) (that \textsc{Severi} calls \emph{%
hyperdifferentiability}) was such a property, the fact recognized by \textsc{Peano}
already in 1884 for the functions of one variable in \cite{peano-lettre,peano-reponse}, and presented later in \cite[(1892)]{peano_mathesis} as an
alternative to usual \emph{differentiability}.

\begin{theorem}[{\textsc{Severi} \cite[(1934)]{Severi-diff}}]
If $A$ is open, then $f\in C^{1}(A)$ if and only if $f$ is strictly
differentiable at every point of  $A$.
\end{theorem}

The next step of \textsc{Severi} was to characterize \emph{strict differentiability}
geometrically for functions with arbitrary (closed) domains. This task was
carried out with the aid of a new concept of tangency, following the same
scheme of geometric characterization of \emph{differentiability}, on
replacing the role of\emph{\ tangent half-lines} by \emph{improper chords}.
\textsc{Bouligand} gave these interrelations in \cite[(1932), pp.\,68-71, 87]{bouligand} (in the
special case where the domain is the Euclidean plane) without furnishing any
precise and complete mathematical formulation\thinspace\footnote{%
\textsc{Bouligand} says in in \cite[(1932), p.\,87]{bouligand}
\par
\begin{quotation}De m\^{e}me que l'hypoth\`{e}se~: r\'{e}duction du contingent \`{a} un plan
pour la surface $z=f(x,y)$, correspond \`{a} la diff\'{e}rentielle prise au
sens de Stolz, de m\^{e}me l'hypoth\`{e}se~: r\'{e}duction du paratingent 
\`{a} un plan pour la surface $z=f(x,y)$, correspond \`{a} la diff\'{e}%
rentielle au sens classique, la fonction $f$ ayant des d\'{e}riv\'{e}es
partielles continues.\end{quotation}

\translation{As the hypothesis of reduction of the contingent to a plane for the surface 
$z=f(x,y)$ corresponds to the [total] differential taken in the sense of
Stolz, the hypothesis of reduction of the paratingent to a plane for the
surface $z=f(x,y)$ corresponds to the differential in the classical sense,
that is, the function $f$ admits continuous partial derivatives.}}.

A linear map $L:\mathbb{R}^{m}\rightarrow \mathbb{R}^{n}$ is a \emph{total
strict differential} of $f$ at an accumulation point $\hat{x}$ of $\limfunc{%
dom}(f)\subset \mathbb{R}^{m}$ provided that $\hat x\in \limfunc{%
dom}(f)$ and %
\begin{equation*}
\lim\nolimits_{x\ne y, x,y\rightarrow \hat{x}}\dfrac{f(y)-f(x)-L(y-x)}{\Vert
y-x\Vert }=0.
\end{equation*}%
\textsc{Severi} provides examples of functions that admit multiple total
differentials and a unique total strict differential\thinspace\footnote{%
For instance \cite[(1944), p.\,283]{Severi-Scorza}, let $A:=\left\{ (x_{1},x_{2})\in 
\mathbb{R}^{2}:\left\vert x_{2}\right\vert \leq x_{1}^{2}\right\} $ and $%
f(x_{1},x_{2}):=0$ for $(x_{1},x_{2})\in A$. Then a total differential $L$
of $f$ at $(0,0)$ fulfills%
\begin{equation*}
\lim\nolimits_{A\ni (x_{1},x_{2})\rightarrow (0,0)}\dfrac{L(x_{1},x_{2})}{%
\Vert (x_{1},x_{2})\Vert }=0,
\end{equation*}%
hence $\left\vert L(x_{1},x_{2})\right\vert \leq \varepsilon \left\vert
x_{1}\right\vert $ for each $\varepsilon >0$, showing that every linear form
such that $L(x_{1},0)=0$ is a total differential. A total hyperdifferential $%
L$ of $f$ at $(0,0)$ satisfies%
\begin{equation*}
\lim\nolimits_{A\ni (y_{1},y_{2}),(x_{1},x_{2})\rightarrow (0,0)}\dfrac{%
L(y_{1}-x_{1},y_{2}-x_{2})}{\Vert (y_{1}-x_{1},y_{2}-x_{2})\Vert }=0.
\end{equation*}%
As for every $\varepsilon >0$ and each $(h_{1},h_{2})$ there exist $%
(y_{1},y_{2}),(x_{1},x_{2})\in A$ and $t>0$ such that $%
(th_{1},th_{2})=(y_{1}-x_{1},y_{2}-x_{2})$, we infer that $\left\vert
L(h_{1},h_{2})\right\vert \leq \varepsilon \left\Vert
(h_{1},h_{2})\right\Vert $, so that $L=0$ is the only total
hyperdifferential of $f$ at $(0,0)$.}. In order to give a geometric
interpretation of total strict differential, \textsc{Severi} makes use of improper
chords, that were also introduced independently by \textsc{Bouligand} \cite[(1928, 1930)]{boul_BSMF,boul_pol}
and called by him \emph{paratingents}. Both \textsc{Severi} and \textsc{Bouligand} consider
the \emph{upper paratangent cone} (\ref{para}) as a family of straight lines
(paratingents,\emph{\ }improper chords). The upper paratangent cone\emph{\ }$%
\limfunc{pTan}\nolimits^{+}(F,\hat x)$ can be characterized in terms of
sequences, as follows: a vector $v\in \limfunc{pTan}\nolimits^{+}(F,\hat x)$
whenever there exist $\left\{ t_{n}\right\} _{n}\rightarrow 0^{+},\left\{
y_{n}\right\} _{n},\left\{ x\right\} _{n}\subset F$ that tend to $\hat x$ such
that%
\begin{equation}
\lim\nolimits_{n}\dfrac{x_{n}-y_{n}}{t_{n}}=v.  \label{paradiff}
\end{equation}

Following \textsc{Guareschi} \cite[(1941), p.\,154]{Guareschi-Roma}, the \emph{linear
paratangent space} of $F$ at $x$ is defined as the linear hull of the
upper paratangent cone of $F$ at $\hat x$.

\begin{theorem}[{\textsc{Severi} \protect\cite[(1934), p.\,189]{Severi-diff}}]
\label{thm:Severi} Let $A\subset \mathbb{R}^{n}$ and $\hat x\in A$ be an
accumulation point of $A$. A function $f:A\rightarrow \mathbb{R}$,
continuous at $\hat x$, is strictly differentiable at $\hat x$ if and only if $%
\limfunc{pTan}\nolimits^{+}(\limfunc{graph}(f),(\hat x,f(\hat x)))$ is included in a
hyperplane without  vertical lines.
\end{theorem}

The \emph{chordal dimension} of \textsc{Guareschi} at an
accumulation point $\hat x$ of a set $F$ is  the dimension of $\limfunc{pTan}\nolimits^{+}(F,\hat x)$.

\begin{theorem}[{\textsc{Guareschi} \protect\cite[(1941), p.\,161]{Guareschi-Roma}}]
\label{thm:guareschi-iper} If the linear paratangent space of $\limfunc{graph}(f)$ at $(\hat x,f(\hat x)))$  does not include  vertical
lines, then there exists a total strict differential if and only if the
chordal dimension of $\limfunc{graph}(f)$ at $(\hat x,f(\hat x))$ is not greater than $%
n$.
\end{theorem}

Analogously to Theorem \ref{thm:sintesi},

\begin{theorem}[{\textsc{Severi} \protect\cite[(1934), p.\,190]{Severi-diff}}]
\label{thm:sintesi-iper} Let $f:A\rightarrow \mathbb{R}^{k}$ where $A\subset 
\mathbb{R}^{m}$, and let $\hat{x}\in A$ be an accumulation point of $A$. Let $L:$ $%
\mathbb{R}^{m}\rightarrow \mathbb{R}^{k}$ be a linear map. Then the
following properties are equivalent:

\begin{enumerate}
\item $f$ is strictly differentiable at $\hat{x}$ and $L$ is a total strict
differential of $f$ at $\hat{x}$;

\item  $f$ is continuous at $\hat{x}$
and $\limfunc{pTan}\nolimits^{+}\left( \limfunc{graph}(f),(\hat{x},f(\hat{x}%
))\right) \subset\limfunc{graph}( L) $;\thinspace\footnote{%
This condition does not appear in \textsc{Severi}, but we evoke it for the sake of
comparison with Theorem \ref{thm:sintesi}.}

\item $\lim\nolimits_{n}\dfrac{f(x_{n})-f(y_{n})}{\Vert x_{n}-y_{n}\Vert }%
=L(v)$ for each $v\in \mathbb{R}^{m}$ and for all sequences $\left\{
x_{n}\right\} _{n},\left\{ y_{n}\right\} _{n}\subset
A$ such that $\lim\nolimits_{n}x_{n}=\hat{x}=\lim\nolimits_{n}y_{n}$, $\lim\nolimits_{n}%
\dfrac{x_{n}-y_{n}}{\Vert x_{n}-y_{n}\Vert }=v$;

\item $L(v)=\lim\nolimits_{\substack{ w\rightarrow v,x\rightarrow \hat{x} 
\\ t\rightarrow 0^{+}}}\dfrac{f(x+tw)-f(x)}{t}$ for every $v\in \limfunc{pTan%
}\nolimits^{+}(A,\hat{x})$.
\end{enumerate}
\end{theorem}

Condition (3) and (4) can be found in \cite[(1934), p.\,190]{Severi-diff} where $L(v)$
fulfilling (3) is called by \textsc{Severi} the \emph{directional hyperderivative} of $f$ at $%
\hat{x}$ along $v$.

\section{Appendix: Turin mathematical community toward Peano}

\label{postcard}

\textsc{Peano}'s interest in logic and in international auxiliary languages coincided
with his progressive marginalization among Turin mathematicians. His
colleagues could not recognize a vital role of \textsc{Peano}'s formal language%
\thinspace\footnote{%
which, among other things, enabled \textsc{Peano} to discover the axiom of choice.}
in the development of mathematics, and were opposed to his teaching methods.
Occurrence of influence groups hostile to \textsc{Peano}'s scientific views led to
his deprivation of the course of calculus, thus of his habitual contacts
with students. Local denigration however did not affect \textsc{Peano}'s worldwide
reputation. He continued to receive highest national distinctions\thinspace\footnote{%
In 1921 the government promoted \textsc{Peano} to \emph{Commendatore} of the Crown of Italy (see \textsc{Kennedy}
\cite[(2006), p.\,215]{kennedy_corrected}).}. Eminent scientists continued to value
him very highly (Appendix \ref{letters}). Nevertheless the persistence of
anti-\textsc{Peano} ambience during his last years, and also for half a century or so
after his death, inescapably left its traces.

\textsc{Tricomi} (1897-1978) joined the faculty of the University of Turin
in 1925. His candidature was strongly supported by \textsc{Peano}'s group and opposed
by the group of \textsc{Corrado Segre} (see \textsc{Tricomi} \cite[(1967), pp.\,18-19]{tricomi}). Here we reproduce a
postcard (and its English translation\thinspace\footnote{%
\translation{Most illustrious professor Giuseppe Peano, of the Royal University of Turin,
4, Barbaroux Street.
\par
Rome, 9th of March 1924
\par
\qquad Illustrious Professor, At the same time that I warmly thank you for
the cordial reception that you wanted to reserve to me [during my visit in
Turin], I have the honour to communicate to you that during the yesterday
meeting of our seminar I spoke to inform the audience about the
conversation, which I was fortunate to have with you on the so called
Zermelo postulate. By the way, I read the passage of your work from the
volume 37 of \emph{Mathematische Annalen} that refers to it, and I had an
impression that all the present were struck by the fact that, eighteen years
before the memoir of Zermelo, you had already formulated, in the very terms
that we use today, the axiom of choice. Moreover Dr Zariski, who studies
here with acuity these things, considered the bibliographical indications
that I got from you, and suggested to relaunch the due revendication of the
contribution of yours and of your school in this difficult area of
mathematics.

\qquad Please accept the finest homages from your devoted F. Tricomi}}) sent
by \textsc{Tricomi} to \textsc{Peano} on the 9th of March 1924.

\begin{quotation}
Illmo Sig$^{\underline{\rm r}}$.

Prof. Giuseppe Peano

della R. Universit\`{a} di Torino

Via Barbaroux, 4

\vspace{0in}------------------------------------

Roma, 9 marzo 1924

\qquad Illustre Professore,

\qquad Nel tempo stesso che vivamente La ringrazio per le cordiali
accoglienze che ha voluto farmi cost\`{a}, mi pregio informarLa che, nella
seduta di ieri del nostro Seminario, ho preso la parola per ragguagliare i
presenti sulla conversazione che ho avuto la fortuna di avere con\ Lei, sul
cos\`{\i} detto postulato di Zermelo.

\qquad Fra l'altro ho letto quel passo del Suo lavoro del t. 37 dei Mathem.\,Annalen che vi si referisce, e mi \`{e} parso che tutti i presenti siano
rimasti colpiti dal fatto che Ella, diciotto anni prima della Memoria di
Zermelo, aveva gi\`{a} formulato, e con le stesse parole che ancora oggid%
\`{\i} usiamo, il principio di scelta.

\qquad Inoltre il Dr. Zarinschi [\emph{sic}], che con acume si occupa qui di queste
cose, ha preso nota delle indicazioni bibliografiche da Lei fornitemi, e si
propone di ritornare su questa doverosa rivendicazione del contributo
portato da Lei e dalla Sua scuola, in questo difficile campo delle
matematiche.

\qquad Voglia gradire, Sig$^{\underline{\rm r}}$ Professore, i pi\`{u} distinti
ossequi del

\qquad \qquad \qquad \qquad \qquad \qquad \qquad Suo devoto F. Tricomi
\end{quotation}

In spite of \textsc{Zariski}'s awareness of \textsc{Peano}'s authorship of the axiom of
choice, we have not found any hint of this fact in the writings of \textsc{Zariski} \cite[(1924-1926)]{zariski_cantor,zariski_zermelo,zariski_zermelo_hilbert,zariski_continuidad,dedekind}. 

\textsc{Tricomi} exercised considerable influence in Turin mathematical
community (and beyond it) till his death. In his writings 
sarcastic  and disdainful opinions on  Italian mathematicians \cite[(1961, 1967)]{tricomiM,tricomi} are profuse. \textsc{Tricomi}
played a decisive role in the discrimination of \textsc{Peano} and used to denigrate
\textsc{Peano} and his school also long after \textsc{Peano}'s death. As reports in \cite%
[pp.\,235-236]{kennedy_corrected} \textsc{Kennedy}, the biographer of \textsc{Peano},

\begin{quotation}
Even later [after 1966] while President of the Academy of Sciences of Turin,
F. G. Tricomi continued to publicly make anti-Peano statements. [...] the
continued attacks on his [Peano] reputation thirty five years later [after
Peano's death] are inexplicable.
\end{quotation}

For a long time the ambiance in Turin (and in Italy) was such that many
preferred to not to reveal their scientific affiliation with the \textsc{Peano}
heritage. Others were simply unaware of the importance of this heritage.

\textsc{Geymonat} (1908-1991), who was graduated in philosophy in 1930 and in
mathematics in 1932 with \textsc{Fubini}, and became an assistant of \textsc{Tricomi}, reports
in \cite[(1986)]{geymonat}:

\begin{quotation}
Quando nel lontano 1934 mi recai a Vienna per approfondire il neopositivismo
di Schlick, portai con me diverse lettere di presentazione (fra le quali
anche una di Guido Fubini); esse vennero accolte favorevolmente e valsero a
creare subito intorno a me una certa cordialit\`{a}. Ma, con mia sopresa, ci%
\`{o} che pes\`{o} pi\`{u} di tutti a mio vantaggio fu il fatto che nel
1930-1931 io ero stato allievo di Peano. Mi sono permesso di ricordare
questo fatto in s\'{e} stesso di nessun rilievo a due scopi: 1) per
sottolineare l'altissima stima di cui Peano godeva, anche dopo la sua morte,
fuori d'Italia; 2) per confessare che purtroppo io pure, come molti altri
giovani appena usciti dall'Universit\`{a} di Torino, non mi rendevo conto
dell'eccezionale valore dell'uomo di cui tuttavia avevo seguito le lezioni
per un intero anno accademico, e col quale avevo avuto tante occasioni per
discorrere anche fuori delle aule accademiche.\thinspace\footnote{%
\translation{When, in the remote 1934, I went to Vienna to study more thoroughly the
neopositivism of Schlick, I carried several recommendation letters (among
which that of Guido Fubini); they were favorably received and created
certain cheerfulness around me. But, to my surprise, what favored me most by
everybody, was the fact that I was Peano's student in 1930-1931. I am
quoting this fact, which is insignificant in itself, for two reasons: 1) to
stress the highest esteem in which Peano was held abroad, also after his
death; 2) to confess that I too, as many other young people graduated from
University of Turin, was not aware of the exceptional worthiness of the man,
the lessons of whom I attended for a whole academic year, and with whom I
had many opportunities to discuss also out of the courses.} 

Presenting himself as a great expert of \textsc{Peano}'s person and works,
\textsc{Geymonat} oscillates between clumsy admiration and commiseration of \textsc{Peano}.}
\end{quotation}

University of Turin has showed little enthusiasm in commemoration of one of
his most illustrious members. \textsc{Kennedy} reports \cite[(2006), p.\,236]%
{kennedy_corrected}:

\begin{quotation}
A few months after his death, the faculty of sciences at the university
considered the possibility of publishing a selection of his writings and
appointed a commission consisting of Carlo Somigliana, Guido Fubini, and F.
G. Tricomi, who worked out a project in 1933. The presence of Tricomi on
this commission practically guaranteed, however, that nothing would come of
the project, and in fact the project was abandoned until after the Second
World War when, Tricomi being in the U.S.A., an analogous project was again
planned by T. Boggio, G. Ascoli, and A. Terracini. In the meantime the
Unione Matematica Italiana [UMI] had decided to publish Peano's work - but
delayed so as not to interfere with the plans of the university. The latter,
however, abandoned this project in 1956 (Tricomi had in the meantime
returned to Turin), so that the UMI then asked Ugo Cassina to propose a
project for publishing Peano's works and on 5 October 1956 named a
commission consisting of Giovanni Sansone, president of the UMI, A.
Terracini, and U. Cassina to make the final selection of works to be
published.
\end{quotation}

The first conference in memory of \textsc{Peano} was organized in 1953 \cite{Cuneo} by \emph{Liceo
Scientifico} of Cuneo, the capital of the province of birth of \textsc{Peano}. 

In 1982 University of Turin organized conference in memory of
\textsc{Peano} for the first time (on the 50th anniversary of \textsc{Peano}'s death).
\textsc{Kennedy}, the biographer of \textsc{Peano}, asked, to no avail, for an invitation \cite[(2006), p.\,IX]{kennedy_corrected}. 
A
booklet of the conference proceedings appeared four years later \cite[(1986)]{celebra}%
. In one of the papers \cite[(1986), p.\,12]{geymonat} of \cite{celebra} \textsc{Geymonat}
recalls the following facts\thinspace\footnote{\translation{In order to save Peano's merits in the area of mathematics, certain persons tried to distinguish two periods [...]. In the first Peano was a talented mathematician, while in the second (decadence phase) his activity was reduced to symbolic logic, passing to linguistic problems related to a search of a universal language [sic] (the pursuit promoted already by Leibniz between seventeenth and eighteenth centuries), the problems that he pretended to able to solve with his latino sine flexione [...]. This was approximately a thesis defended by Fubini, his great antagonist at the Faculty of Turin, during a talk held at the mathematical Seminar of this faculty about 1930, I do not remember exactly, but in any case when Peano was still alive. But even that talk did not succeed to reconcile the positions of Fubini and Peano [...].}}:

\begin{quotation}
Per poter salvare i meriti di Peano nel campo matematico, alcuni avevano
cercato di distinguere nettamente due fasi [...]. Nella prima fase Peano
sarebbe stato un valente matematico, mentre nella seconda (o fase della
decadenza) egli si sarebbe ridotto a occuparsi di logica simbolica, passando
poi a problemi linguistici connessi alla ricerca di un linguaggio universale
[\emph{sic}] (ricerca gi\`{a} promossa da Leibniz negli anni a cavallo fra il Sei e
il Settecento), problemi che egli ritenne di poter risolvere con il suo
latino sine flexione [...]. Questa all'incirca fu la tesi sostenuta da
Fubini, il suo grande avversario nella Facolt\`{a} di Torino, in una
conferenza tenuta al Semi\-nario matematico di tale Facolt\`{a}, non ricordo pi%
\`{u} esattamente se poco prima o poco dopo il 1930, comunque mentre Peano
era ancora in vita. Ma neanche questa conferenza riusc\`{\i} a conciliare le
due posizioni di Fubini e Peano [...].
\end{quotation}

Recalling events of that conference in \cite[(1982)]{lolli}, \textsc{Lolli}, who
graduated with \textsc{Tricomi} in 1965 and became an assistant of \textsc{Geymonat} in 1967, alludes to a \emph{curtain of silence} of
the Turin mathematical community around the \emph{embarrassing and bizarre
personage who, for about fifty years, disturbed and discomfitted, and in the
last thirty years almost dishonored the whole profession}\thinspace\footnote{%
\textsc{Lolli}'s words:
\quote{la cortina di silenzio [of the Turin mathematical cummunity around
the] [\dots] scomodo e bizzarro personaggio che per circa cinquanta anni aveva
disturbato ed imbarazzato, e negli ultimi trenta quasi disonorato la intera
professione.}
}. In his book \cite[(1985), p.\,8]{lolli2}, \textsc{Lolli} qualifies
\textsc{Peano} as a \emph{pathetic inventor of symbols} and, in the same book \cite[(1985),
p.\,50]{lolli2},\emph{\ who made through cowardice the great refusal}%
\thinspace\footnote{\textsc{Dante} \cite[Inferno, Canto III]{dante}: ``Colui che fece per viltade il gran rifiuto''.} in reference to \textsc{Dante}'s \emph{Divina
Commedia}.  \footnote{Ironically, in 2000 \textsc{Lolli} was recipient  of  a Peano Prize,  sponsored by Department of Mathematics of Turin.}

The persistence of anti-Peano ambience in Turin Mathematical Community  a half century after \textsc{Peano}'s death, was nourished and reinforced by a surprisingly poor knowledge of his works. In \cite[1959]{geymonat59} \textsc{Geymonat}, an authoritative member of that community, on the occasion of the edition of \textsc{Peano}'s \emph{Selected Works}  by \textsc{Cassina}, wrote\thinspace\footnote{%
\translation{The second volume [of Peano's Selected Works]  [\dots] gathers works in mathematical logic [\dots] [and] in interlingua and algebra of grammar. This juxtaposition [\dots] confirms without doubt \textsc{Cassina}'s opinion, after which mathematical logic and linguistic research constitute, in \textsc{Peano}, two phases [\dots] of the same grand program designed to realize [\dots] the teaching of \textsc{Leibniz}.

This thesis is of particular importance, because it undermines the legend [\emph{sic}], following which the linguistic interests of \textsc{Peano} would be a fruit of his senile decadence.}}:

\begin{quotation} Il II volume [delle Opere Scelte di Peano] [\dots] raccoglie lavori di logica matematica [\dots] [e] lavori di interlingua ed algebra della grammatica. L'accostamento [\dots] conferma in modo incontestabile l'opinione di \textsc{Cassina}, secondo cui logica matematica e ricerche linguistiche costituiscono, in \textsc{Peano}, due fasi [\dots] di un medesimo grandioso programma volto a realizzare [\dots] l'insegnamento leibniziano. 

La tesi ha una particolare importanza, perch\'e sfata la leggenda [\emph{sic}] secondo cui gli interessi linguistici peaniani sarebbero stati il frutto di una decadenza senile del Nostro.
\end{quotation}

 Multiple contributions of \textsc{Mangione} on the history of logic to the six volumes of \textsc{Geymonat}'s \emph{Storia del pensiero filosofico e scientifico} \cite[(1971-1973)]{geymonat-storia} indicate persisting poor knowledge of \textsc{Peano}'s works. \textsc{Mangione}'s contributions, very much appraised by Italian logicians and philosophers, are completely unknown to mathematicians. They were collected in \emph{Storia della logica} \cite[(1993)]{mangione} a few years ago, without any change of attitude with regard to \textsc{Peano} and his School, who are ridiculed therein.

In  \emph{La Stampa}, a daily of Turin, in October 1995 \textsc{R.\,Spiegler} declared that certainly \textsc{Peano} spent some periods in a madhouse. This news without any basis was belied by \textsc{Lalla Romano}, a \textsc{Peano}'s great-niece. A mathematician and our colleague asked \textsc{Spiegler} (who is also a mathematician) where he took this absurd information;   \textsc{Spiegler} replied that he had learned this from \textsc{G.-C.Rota} who, in turn, was informed by nobody else but \textsc{Tricomi} in person. \footnote{\textsc{Rota} wrote in \emph{Indiscrete Thoughts} (Birk\"auser, 1997, page 4): ``Several outstanding logicians of the twentieth century found shelter in asylums at some time in their lives: Cantor, Zermelo, G\"odel, Peano, and Post are some.''

Another example of a disdainful attitude toward \textsc{Peano} was the adjectival use of ``peanist'' rather than of more standard and graceful ``peanian''.  The word ``peanist'' was introduced by the renowned  historian  \textsc{Grattan-Guinness}; it evokes the word  ``opportunist'' that was used in a judgement of  \textsc{Grattan-Guinness} on \textsc{Peano}'s works: ``Both in his mathematics and his logic, he [\textsc{Peano}] seems to me to have been an opportunist'' \cite[(1986)]{opportunist}.}

More recently University of Turin edited \emph{Opera omnia} \cite[(2002)]{omnia};
\textsc{Peano} is the celebrity whom \emph{Accademia delle Scienze} of Turin put on its home
page

\begin{center}
\texttt{http://www.torinoscienza.it/accademia/home}.
\end{center}

An international congress \emph{Giuseppe Peano e la sua Scuola, fra
matematica, logica e interlingua}\ commemorating the 150th anniversary of
\textsc{Peano}'s birth and 100th anniversary of Formulario Mathematico took place in
Turin in October 2008 at the Academy of Science of Turin and the Archive of
State.

\textsc{Peano}'s is not the first case of an ostracism against a mathematical
precursor. As in other cases, the resulting prejudice is inestimable. And,
as a rule, pupils cannot expect a better destiny.

A famous economist \textsc{Luigi Einaudi} (1874-1961), who was a professor of University of Turin
before becoming the president\thinspace\footnote{%
from 1948 to1955.} of the Italian Republic, witnesses in 1958 \cite{einaudi}:

\begin{quotation}
Il professor Peano fu vero maestro, sia per l'invenzione di teoremi, che
ritrovati poi da altri, resero famosi gli scopritori, sia per l'universalit%
\`{a} del suo genio. Nemmeno a farlo apposta, taluni suoi assistenti ai
quali si pronosticava un grande avvenire nel campo matematico, presero
tutt'altra via. [\dots] Vacca [assistente di Peano], divenne [\dots] professore
universitario di lingua e letteratura cinese [\dots]. [Vailati] nonostante la
crescente estimazione in cui era tenuto nel mondo scientifico italiano e
straniero, [\dots] non ottenne la cattedra alla quale doveva aspirare. [\dots]
Cos\`{\i} fu che Vailati scomparve dall'orizzonte torinese per girare
l'Italia come insegnante nelle scuole medie.\thinspace\footnote{%
\translation{Professor Peano was a real master, as for the invention of theorems, which
rediscovered later by others, made them famous, as for the universal
character of his genius. Not deliberately, several of his assistants, who
had great prospects in mathematics, took completely different ways. [\dots]
Vacca, [an assistant of Peano] became [\dots] a university professor of Chinese
language and literature [\dots]. [Vailati] who despite the growing esteem in
which he was held by Italian and foreign scientists [\dots] did not obtain a
professorship, for which he could legitimately pretend. [\dots] So Vailati
disappeared from the Turin horizon to move around Italy as a secondary
school teacher.}}
\end{quotation}

\section{Appendix: International mathematical community toward Peano}

\label{letters}

Despite the depicted ambience at the University of Turin, \textsc{Peano} was held in
high esteem by numerous famous scientists also in that period.\thinspace\footnote{%
A writer \textsc{Lalla Romano} (1906-2001), \textsc{Peano}'s great-niece describes the atmosphere of
\textsc{Peano}'s house, where she was a guest (1924-1928) during her
unversity studies \cite[(1979), p.\,8]{lalla}:
\par
[...] lo zio [Peano] riceveva le visite: studenti, per lo pi\`{u} stranieri
- perfino cinesi - ossequiosissimi, dal sorriso esitante, l'inchino a
scatto; e scienziati [...] guardavano lo zio con venerazione. Mentre lui,
cupo, la barba arruffata, andava avanti e indietro nella stanza, scuotevano
la testa.
\par
\translation{[...] my uncle [Giuseppe Peano] received visitors: students, mostly
foreigners -- even Chinese -- obsequious, smiling hesitatingly, bowing
snappingly; scientists [...] looked at my uncle with veneration. While he,
gloomy, with his ruffled beard, walked to and fro, they shaked their heads.}}

Among the letters and telegrams sent to \textsc{Peano} on his 70th birthday are those
of \textsc{Guareschi}, \textsc{Dickstein}, \textsc{Zaremba}, \textsc{Fr\'{e}%
chet}, \textsc{Hadamard}, \textsc{Tonelli}  and   \textsc{Levi-Civita}   \cite[(1928)]{inter}]. 

We include few samples  of letters
and other signs of recognition around 1930. They are extracted from a \cite[(2002)]{archivio}.

\begin{center}------------\end{center}

A letter from \textbf{Benjamin Abram Bernstein} (1881-1964)

\begin{quotation}
University of California, Department of Mathematics, Berkeley, California,
Feb. 8, 1928

My dear Professor Peano -

I am anxious to get the Rivista di Matematica v. 1-8, and the Formulaire Math%
\'{e}matique, v. 1-5. I shall appreciate it greatly if you can tell me if
these can be still got from the publishers and at what price.

\qquad With keen appreciation of your great work in logic, I am,

\qquad Sincerely yours, BABernstein.
\end{quotation}

\begin{center}------------\end{center}

A letter\thinspace\footnote{%
\translation{Warsaw, 31.VII.1928
\par
Dear Professor, Please forgive me that I write in German, but unfortunately
I do not know that much Italian in order to communicate with you in your
mother tong.

I did not expect at all that I would be able to take part in the
International Congress of Mathematicians in Bologna. Only now I have this
possibility. Therefore, I ask you, if it were still in some way possible to
accept a delayed registration of my communications. For years I have been working in
the area of mathematical logic, but I have not yet published my most
important results on the propositional calculus and its history. I would be
delighted if I could present my results in Italy, that has so many merits in mathematical logic,  to the international
learned audience.
\par
If it were no longer possible that I actively participate in the congress, I
would be very grateful for information about it.
\par
Please, accept the expression of my greatest respect.
\par
Dr.\'Jan \L ukasiewicz, Professor Philosophy and a former Rector of Warsaw
University /Poland/.
Address: Prof. Dr. J.\ \L ukasiewicz, Brzozowa street, 12, Warsaw, Poland.}}
from \textbf{Jan \L ukasiewicz} (1878-1956)

\begin{quotation}
Warszawa, 31.VII.1928

Sehr Geehrter Herr Professor!

\qquad Bitte mich vielmals zu entschuldigen, dass ich deutsch schreibe, aber
ich verstehe leider nicht soviel italienisch, um mich mit Ihnen in Ihrer
Mutterssprache zu verst\"{a}ndigen.

\qquad Ich habe gar nicht gehofft, dass ich an dem Internationalen Kongresse
der Mathematiker in Bologna werde teilnehmen k\"{o}nnen. Nun hat sich mir
die M\"{o}glichkeit geboten, nach Bologna zu kommen. Ich bitter daher, Herr
Professor, wenn es nur irgendwie m\"{o}glich ist, meine versp\"{a}tete
Anmeldung von Kommunikaten g\"{u}tigst ber\"{u}cksichtigen zu wollen. Seit
Jahren arbeite ich im Gebiete der mathematischen Logik, doch habe ich meine
wichtigsten Ergebnisse aus dem Aussagenkalk\"{u}l und dessen Geschichte
bisher nicht ver\"{o}ffentlicht. Es w\"{a}re mir sehr lieb, wenn ich meine
Resultate gerade in Italien, das so sehr f\"{u}r die mathematische Logik
verdient ist, der internationalen Gelehrtenwelt vorlegen k\"{o}nnte.

\qquad Sollte es nicht mehr m\"{o}glich sein, dass ich am Kongresse aktiv
teilnehme, so w\"{a}re ich f\"{u}r eine Mitteilung dar\"{u}ber sehr dankbar.

\qquad Bitte, Herr Professor, den Ausdruck meiner vorz\"{u}glichsten
Hochachtung entgegenzunehmen

Dr. Jan \L ukasiewicz, Professor f\"{u}r Philosophie und gewesener Rektor
der Universit\"{a}t Warschau /Polen$/.$

Adresse: Prof. Dr. J.\ \L ukasiewicz, Warszawa, Brzozowa 12. /Varsovia
[\emph{sic}], Polonia/
\end{quotation}

\begin{center}------------\end{center}

In a speech at the Congress of Mathematicians in Bologna on the 3rd of
September 1928 \cite[p.\,4]{hilbert} \textbf{David Hilbert} (1862-1943) talks about
\textsc{Peano}'s symbolic language\thinspace\footnote{%
\textsc{Peano} did not participate in that Congress because of his brother's death.}:

\begin{quotation}
[...] ein wesentliches Hilfsmittel f\"{u}r meine Beweistheorie [ist] die
Begriffsschrift; wir verdanken dem Klassiker dieser Begriffsschrift, Peano,
die sorgf\"{a}ltigste Pflege und weitgehendste Ausbildung derselben. Die
Form in der ich die Begriffsschrif brauche, ist wesentlich diejenige, die
Russell zuerst eingef\"{u}rht hat.\thinspace\footnote{%
\translation{[...] an essential tool for my proof theory is ideography; we owe to the
classical author of this ideography, Peano, most thorough care and utmost
cultivation of it. The form, in which I use this ideography, is essentially
that Russell has first introduced.}}
\end{quotation}

\begin{center}------------\end{center}

A letter\thinspace\footnote{%
\translation{Pisa, 12 January 1931=IX%
${{}^\circ}$%
, Illustrious Professor,
\par
During this year 1931 the \textquotedblleft Annals of the Scuola Normale
Superiore\textquotedblright\ of Pisa will absorb the Annals of Tuscan
universities and will be transformed in a great international periodical, of
the type of \textquotedblleft Annals of the \'{E}cole Normale Sup\'{e}%
rieure\textquotedblright\ of Paris. The mathematical section, that will
receive memoirs and notes of excellent Italian and foreign scientists, will
appear, each year, in four volumes of 100 pages each.
\par
\qquad The Scuola Normale hopes to count you among the collaborators of the
so renewed Annals; and I particularly would be very glad if I could include
one of your papers in the first volumes of the new series.
\par
\qquad Would you be so kind to gratify me? With anticipated thanks and many
homages. Your most devoted, L.\,Tonelli}} from \textbf{Leonida Tonelli} (1885-1946)

\begin{quotation}
Pisa, 12 gennajo [sic] 1931=IX%
${{}^\circ}$%

Illustre Professore,

Nel corrente anno 1931, gli \textquotedblleft Annali della Scuola Normale
Superiore\textquotedblright\ di Pisa assorbiranno gli \textquotedblleft
Annali delle Universit\`{a} Toscane\textquotedblright\ e si trasformeranno
in un grande periodico internazionale, del tipo degli \textquotedblleft
Annales de l'\'{E}cole Normale Sup\'{e}rieure\textquotedblright\ di Parigi.
La parte matematica, che accoglier\`{a} Memorie e Note di valorosi
scienziati italiani e stranieri, si presenter\`{a}, ogni anno, con quattro
fascicoli, ciascuno di 100 pagine.

\qquad La Scuola Normale spera di poterLa annoverare fra i collaboratori
degli Annali cos\`{\i} rinnovati; ed io, in particolare, sarei molto lieto
se potessi inserire un Suo lavoro nei primi fascicoli della nuova serie.

\qquad Vuole essere tanto gentile da accontentarmi?

\qquad Con anticipati ringraziamenti e molti ossequi.

\qquad Suo devotissimo, L. Tonelli
\end{quotation}

\begin{center}------------\end{center}

A letter\thinspace\footnote{%
\translation{Warsaw, 2nd of November 1932. Dear Professor, I take freedom to bother you
with my personal affairs. I have namely a prospect, for the coming year
1933/4, to obtain the Rockefeller fellowship to study abroad, and would be
very glad if I could work sometime in Turin under your supervision. Would
you kindly agree to this?
\par
Looking forward to your kind reply, I remain in deep respect.
\par
Dr.\,A. Tarski, Private Docent at Warsaw University (Poland, Warsaw XXI, Su\l %
kowskiego street 2 app.5)}} from \textbf{Alfred Tarski} (1902-1983)

\begin{quotation}
Warschau, 2.XI.32\thinspace\footnote{%
Peano died on the 20th of April 1932.}

\qquad \qquad Hoch verehrter Herr Professor!

\qquad Ich nehme mir die Freiheit, Sie mit einer privaten Angelegenheit zu
behelligen. Ich habe n\"{a}mlich die Aussicht, f\"{u}r das kommende Jahr
1933/4 das Rockefeller-Stipendium f\"{u}r das Studium in Ausland zu
bekommen, und w\"{u}rde mich sehr freuen, wenn ich eine Zeit unter Ihrer F%
\"{u}hrung in Turin arbeiten durfte. W\"{u}rden Sie damit einverstanden sein?

\qquad In Erwartung Ihrer freundlichen Antwort verbleibe ich inzwischen in
vorz\"{u}glicher Hochachtung

\qquad \qquad Dr. A. Tarski, Privat-Dozent a.d. Universit\"{a}t Warschau
(Polonia, Warszawa XXI, ul. Su\l kowskiego 2 m.5)
\end{quotation}

\section{Appendix: Biography of Giacinto Guareschi}

\label{biographies}

\emph{We provide a somewhat detailed biography of Guareschi, because it is
not available, except a brief mention in Atti dell'Accademia Ligure} \cite%
{rizzi}.\emph{\ Biographies of other mathematicians referred to in this
paper are easily obtainable.}

\textbf{Giacinto Guareschi} (1882-1976) was born in Turin on the 2nd of
October 1882. His father, \textsc{Icilio} (1847-1918) was a famous
chemist-pharmacologist \cite{icilio}, a member of Accademia delle Scienze of Turin at the
same time as Peano. His mother was \textsc{Anna Maria Pigorini} (\dag\ 1942).
\textsc{Guareschi} had a sister \textsc{Paolina} and a brother \textsc{Pietro} (1888-1965), a
distinguished chemical engineer, member of \emph{Accademia Ligure}.

\textsc{Guareschi} studies mathematics at the University of Turin graduating in 1904.
In a letter of 1932 \cite[p.\,87]{vacca} to \textsc{Vacca}, he recalls the importance of \textsc{Severi} and \textsc{Vacca}
(assistants of, respectively, \textsc{D'Ovidio} and \textsc{Peano}) for his mathematical
education. He was assistant of projective geometry at the University of
Turin (1904-1906), and of analytic geometry at the University of Pavia
(1907-1910)\thinspace\footnote{%
at a suggestion of  \textsc{Berzolari}.}. In 1910 he obtained a professorship of
high school \emph{(liceo)} to voluntarily retire in 1944 in order not to
collaborate with, and to avoid to swear faithfulness to the Fascist regime.
During his high school teacher carrier, \textsc{Guareschi} served as a principal and
was appointed\thinspace\footnote{%
by the minister of National Education, without having asked for it.
\textsc{Guareschi} was not happy with this nomination, mainly because it interfered
with his research (namely, on differentiability and tangency), but could not
refuse due to the legal system at that moment. Soon after he realized that
the \textsc{Mussolini} government politicized education. In \textsc{Gareschi}'s words:
\par
\begin{quotation}
[il] pagliaccio di Predappio [aveva reso la carica di Provveditore]
squisitamente politica
\end{quotation}
\par
\translation{because the clown of Predappio made this position exquisitly political [The
reference to Mussolini who was born in Predappio]}.
\par
Contrary to \textsc{Guareschi}, \textsc{Severi} is an enthousiastic follower of \textsc{Mussolini} (see \textsc{Guerraggio}-\textsc{Nastasi} in \cite[(1993, 2005)]{guerra_camicia,guerra_gentile}.} a \emph{provveditore }\thinspace\footnote{%
a provincial responsible of education.} in July 1936. From November 1936
\textsc{Guareschi} continued to ask to be exempted\thinspace\footnote{%
The reason was primarily political, because \textsc{Guareschi} was opposed to the
Fascit regime, however he could not openly evoke it, as this would amount to
severe persecution.}, and, after several refusals, was finally dismissed in
1938.

On the 21st of November 1914 he was enrolled in the army and participated in
the First World War. He left the army on 15th of May 1919 with the grade of
captain; in 1921 he was granted a commemorative medal of the First World
War. In 1931 he was promoted to the grade of major of artillery, and on 11th
of June 1940 was enrolled to the army to be demobilized on the 19th of
August of the same year with the grade of lieutenant-colonel.

In 1924 \textsc{Guareschi} started pedagogical activity in projective and analytic
geometry at the University of Genoa, where he became a \emph{libero docente }%
\thinspace\footnote{%
The title of \emph{libero docente}, granted on the basis of scientific
publication, entitled to teach courses at a university.} of algebra on 13th
of March 1929. He kept this position till 1952 (when he became 70, which was
the legal retirement age). Due to a derogation, he taught at the University
of Genoa till 1959.

In 1927 \textsc{Guareschi} was elected a corresponding member of \emph{Accademia} Ligure di
Scienze e Lettere (proposed by \textsc{Loria} and \textsc{Severini}) and in 1957
its effective member. In 1956 he and his brother \textsc{Pietro} donated to \emph{Accademia}
manuscripts of their father \textsc{Icilio}.

\textsc{Guareschi} married \textsc{Gemma Venezian} (1897-1975). Their only son, \textsc{Marco}, was
born on the 21st of March 1922. In 1944 he joined the underground army,
which was, in terms used by \textsc{Guareschi}, \emph{la sola via dell'onore}
(the only way of honor). On the 11th of April 1944 \textsc{Marco} was arrested%
\thinspace\footnote{%
at the \emph{rastrellamento }(sweep) of \emph{Benedicta}, where more than
hundred partisans were executed and other 400 arrested. \textsc{Guareschi}
reconstructed the event in \cite[(1951)]{guareschiX3}, which became a basic source
for \cite[(1967)]{pansa} of \textsc{Pansa}.} and deported to Germany where he died
in a concentration camp in April 1945\thinspace\footnote{%
First to Mauthausen, later in August 1944 to Peggau (near Graz) and finally
to the so called Russian Camp where he died between 10 and 12 April 1945.}.
The pain of \textsc{Guareschi} and his wife was amplified by uncertainty about their
son's fate, as, for a couple of years, they did not have reliable
information about his passing. Since then \textsc{Guareschi} dedicated himself to
promotion to reconstruction of the history of the \emph{Resistenza}
(Italian underground army) and to defence of its values; in doing so, he
collaborated with several Italian and international associations\thinspace\footnote{%
For example, \emph{Istituto storico della Resistenza} in Liguria, ANED
(\emph{Associazione nazionale ex deportati}), ANPI (\emph{Associazione nazionale
partigiani d'Italia}), ANCR (\emph{Associazione combattenti e reduci}), ANPPIA
(\emph{Associazione Nazionale Perseguitati Politici Italiani Antifascisti}),
\emph{Consiglio Federativo della Resistenza}, \emph{Conseil Mondial de la Paix}.}.

The postwar years were extremely difficult for \textsc{Guareschi} and his wife.
\textsc{Guareschi} had neither salary nor pension, because he resigned from the
public service during the Fascist period. In January 1946 \textsc{Guareschi} wrote

\begin{quotation}
Io me ne sono andato [dalla scuola] per non servirla [la repubblica
fascista] al tempo dell'obbligo del giuramento, e nemmeno ho giurato agli
Ufficiali in congedo; n\'{e} pi\`{u} ho esercitato l'incarico Universitario,
sfidando la fame. [...] Sono agli estremi dal lato finanziario; i mesi
arretrati [per il pagamento dello stipendio e della pensione] sono ormai 21.%
\thinspace\footnote{%
\translation{I quit [the teaching] in order not to serve [the fascit republic] at the
time of obligation of oath of allegiance, nor I swore as an officer in
leave; nor I had a university appointment, defying the hunger. [...]
Financially I am destitute. The arrears [of wage and pension] are already
for 21 months.}}
\end{quotation}

\textsc{Guareschi} successfully applied to be readmitted as a high school professor,
because the political nature of his resignation in 1944 was recognized.

After the war \textsc{Guareschi} had various political commitments. In 1945 he became
a mayor of a village Serravalle Scrivia (Alessandria). In 1953 he was an
unsuccessful candidate (from the lists of PCI\thinspace\footnote{%
Italian Communist Party.}) for senator. In recognition of their intense
political activity, \textsc{Giacinto} and \textsc{Gemma} \textsc{Guareschi} received a gold medal in
1956. On his retirement from the secondary education, on the 28th of
September 1950, three principal newspapers of Genoa (\emph{Il lavoro nuovo}, 
\emph{Il secolo XIX} and \emph{l'Unit\`{a}}) published a paper about
\textsc{Guareschi}, writing, among other things,

\begin{quotation}
Inflessibile nei riguardi delle ingerenze del regime fascista nella vita
della scuola, durante la lotta contro i nazifascisti ha offerto alla Patria
l'unico figlio barbaramente trucidato a Mauthausen.\thinspace\footnote{%
\translation{[Guareschi] inflexible with respect to the intrusions in the school life of
the Fascist regime, during the fight against Nazifascists, [he] offered to
the Fatherland his only son, barbarically slain in Mauthausen.}}
\end{quotation}

\textsc{Giacinto Guareschi} died on the 9th of August 1976 in Serravalle Scrivia near
Alessandria, in a poor country house, where he lived his last years. Various
scholarships, prizes were founded and monuments were erected in memory of
\textsc{Guareschi}.

Mathematical interests of \textsc{Guareschi} are principally geometry and algebra,
and starting from 1934, differentiability and tangency  (see previous Sections \ref{diff-cara} and \ref{strictdiff-cara}) and, finally, 
characterization of smooth manifolds (see \textsc{Greco} \cite{manifold} for details). \textsc{Guareschi}'s works are reviewed in  \emph{Jahrbuch \"uber die Fortschritte de Mathematik} (JFM), in \emph{Zentralblatt Math} (Zbl)  and in \emph{Mathematical Reviews}(MR).\thinspace\footnote{The reviewers of papers on differentiability and tangency are \textsc{H.\,Kneser}, \textsc{O.\,Haupt}, \textsc{A.B.\,Brown}, \textsc{H.\,Busemann}, \textsc{G.\,Scorza Dragoni}, \textsc{O.\,Zariski}, \textsc{T.\,Viola} and \textsc{A. Gonz\'alez Dom\'\i nguez}.}

Scientific publications of \textsc{Guareschi} cease with the death of his son. Nevertheles his interest for mathematics persists during all his life.
In his nineties he collaborates with \textsc{G.\,Rizzitelli} on the edition of a collection of applications of mathematics,  and announces to the secretary of \emph{Accademia Ligure} his intention to publish a paper on algebra.  
\textsc{Guareschi} wrote 3 books for didactic use \cite{rizzi} and  35 mathematical papers. The following bibliography contains only mathematical papers and 5 writings on the \emph{Resistenza}.

\renewcommand{\refname}{{\textsc{Guareschi's Bibliography}}}

\section{Appendix: A chronological list of mathematicians}

For reader's convenience, we provide a chronological list of some
mathematicians mentioned in the paper, together with biographical sources.

The \texttt{html} file  with biographies of mathematicians listed below with an asterisk can be attained at University of St
Andrews's web-page 

\texttt{http://www-history.mcs.st-and.ac.uk/history/\{Name\}.html}
\medskip

\textsc{Descartes}, Ren\'e (1596-1650) (*)

\textsc{Huygens}, Christiaan (1629-1695) (*) 

\textsc{Leibniz}, Gottfried Wilhelm (1646-1716) (*)

\textsc{Grassmann}, Hermann (1809-1877) (*)

\textsc{Genocchi}, Angelo (1817-1889) (*) 

\textsc{Serret}, Joseph (1819-1885) (*)

\textsc{Bertrand}, Joseph L.F. (1822-1900) (*) 

\textsc{Jordan}, Camille (1838-1922)  (*) 

\textsc{Thomae}, Carl J. (1840-1921) (*) 

\textsc{Stolz}, Otto (1842-1905)  (*)

\textsc{D'Ovidio}, Enrico (1842-1933), see \textsc{Kennedy} \cite{kennedy_corrected}

\textsc{Schwarz}, Hermann A. (1843 - 1921) (*) 

\textsc{Dini}, Ulisse (1845-1918) (*)

\textsc{Klein}, Felix (1849-1925) (*)

\textsc{Dickstein}, Samuel (1851-1939) (*)

\textsc{Peano}, Giuseppe (1858-1932) (*), see  \textsc{Kennedy} \cite{kennedy_corrected} 

\textsc{Hilbert}, David (1862-1943) (*)

\textsc{Loria}, Gino B. (1862-1954) (*)

\textsc{Segre}, Corrado (1863-1924) (*)

\textsc{Hadamard}, Jacques S. (1865-1963) (*)

\textsc{Saks}, Stanislaw (1897-1942) (*)

\textsc{Couturat}, Louis (1868-1914) (*)

\textsc{Hausdorf}, Felix (1868-1942) (*)

\textsc{Zermelo}, Ernst (1871-1953) (*)

\textsc{Severini}, Carlo (1872-1951),   see \emph{Boll. Un. Mat. It.} \textbf{7}:98--101, 1952

\textsc{Russell}, Bertrand  (1872-1970) (*)

\textsc{Levi-Civita}, Tullio  (1873-1941) (*)

\textsc{Levi}, Beppo (1875-1961),  see  \textsc{Kennedy} \cite{kennedy_corrected}

\textsc{Vacca}, Giovanni (1875-1953) (*) 

\textsc{Boggio}, Tommaso (1877-1963) (*) 

\textsc{Fr\'{e}chet}, Maurice (1878-1973) (*)

\textsc{\L ukasiewicz}, Jan (1878-1956) (*)

\textsc{Fubini}, Guido (1879-1943) (*)

\textsc{Severi}, Francesco (1879-1961) (*) 

\textsc{Bernstein}, Benjamin A. (1881-1964),  see \emph{Univ.\,California:\,In Memoriam, 1965}

\textsc{Guareschi}, Giacinto (1882-1976)

\textsc{Tonelli}, Leonida (1885-1946),  see \textsc{Tonelli},\,\emph{Opere Scelte}, Cremonese, 1963

\textsc{Bouligand}, Georges (1889-1979), see \texttt{http://catalogue.bnf.fr}

\textsc{Wilkosz}, Wiltold (1891-1941), see 
\texttt{http://www.wiw.pl/matematyka/Biogramy}

\textsc{Ackermann}, Wilhelm (1896-1962) (*)

\textsc{Cassina}, Ugo (1897-1964), see  \textsc{Kennedy} \cite{kennedy_corrected}

\textsc{Tricomi}, Francesco G. (1897-1978) (*)

\textsc{Saks}, Stanislaw (1897-1942) (*)

\textsc{Zariski}, Oscar (1899-1986) (*)

\textsc{Tarski}, Alfred (1902-1983) (*)

\textsc{Segre}, Beniamino (1903-1977) (*)

\textsc{Geymonat}, Ludovico (1908-1991),  see \texttt{www.torinoscienza.it/accademia/}

\textsc{Choquet}, Gustave (1915-2006), see \emph{Gazette des Math.} v111:74-76, 2007.

\textsc{Mangione}, Corrado (1930-2009), see \texttt{http://dipartimento.filosofia.unimi.it/}

\textsc{Rota}, Gian-Carlo (1932-1999) (*)

\renewcommand{\refname}{{\textsc{References}}}


\end{document}